\documentclass[11pt]{article}
\usepackage{enumerate}
\usepackage{epsfig,alltt}
\usepackage[english]{babel}
\usepackage{blindtext}
\usepackage{dsfont}
\usepackage{amssymb}
\usepackage{amsfonts}
\usepackage{graphicx}
\usepackage{amsmath}
\usepackage{amsthm}
\usepackage[colorlinks=true,linkcolor=blue,citecolor=blue,pdfborder={0 0 0}]{hyperref}
\usepackage{natbib}
\usepackage[utf8]{inputenc}
\usepackage{amsthm,amsmath}
\usepackage[left=2.5cm,right=2.5cm,top=2.5cm,bottom=2.5cm]{geometry}
\usepackage{graphicx}
\usepackage{subcaption} 

\usepackage{booktabs}
\usepackage{multirow}
\usepackage{siunitx}
\usepackage{floatflt}
\usepackage{makeidx}
\usepackage{algorithm}

\usepackage{algorithmic}

\usepackage{url}
\usepackage{float} 
\theoremstyle{remark}

\usepackage{soul}

\usepackage{color}
\usepackage{appendix}

\newcommand{\begeq}[1]{\begin{equation} \label{#1}}
	\newcommand{\fineq}{\end{equation}}

\newcommand{\E}{\mathrm{E}}

\newcommand{\Var}{\mathrm{Var}}

\newcommand{\Cov}{\mathrm{Cov}}

\title{Hybrid estimation for a mixed fractional Black-Scholes model with random effects from discrete time observations}

\author{Nesrine CHEBLI$^{1,2}$\footnote{corresponding author}
	\footnote{e-mail address: nesrine.chebli@univ-poitiers.fr}, Hamdi FATHALLAH$^{2}$\footnote{e-mail address: hamdi.fathallah@essths.u-sousse.tn} and Yousri SLAOUI$^{1}$\footnote{e-mail address: Yousri.Slaoui@math.univ-poitiers.fr} \\
	$^{1}$ Laboratory of Mathematics and Applications, University of Poitiers, France\\
	$^{2}$ Laboratory of Mathematics Deterministic and Random Modeling, University of Sousse,\\ Tunisia
}

\numberwithin{equation}{section}
\begin{document}
\newtheorem{theorem}{\bf Theorem}[section]
\newtheorem{cor}{\bf Corollary}[section]
\newtheorem{theor}{\bf Theorem}
\newtheorem*{theo}{\bf Theorem}
\newtheorem{prop}{ \bf Proposition}
\newtheorem{lem}{\bf Lemma}
\newtheorem{coro}{\bf Corollary}
\newtheorem{ex}{\bf Example}
\newtheorem{prof}{\bf Proof}
\newtheorem{defi}{\bf Definition}
\newtheorem{rem}{\bf Remark}
\newtheorem*{nota}{\bf Notations}
\newtheorem*{assump}{\bf Assumptions}
\newtheorem*{com}{\bf Comments on the assumptions}
\newtheorem*{app}{\bf Appendix}
\date{ }
\maketitle

\textit{\bf Abstract}
We propose a hybrid estimation procedure to estimate global fixed parameters and subject-specific random effects in a mixed fractional Black-Scholes model based on discrete-time observations. Specifically, we consider $N$ independent stochastic processes, each driven by a linear combination of standard Brownian motion and an independent fractional Brownian motion, and governed by a drift term that depends on an unobserved random effect with unknown distribution. Based on $n$ discrete time statistics of process increments, we construct parametric estimators for the Brownian motion volatility, the scaling parameter for the fractional Brownian motion, and the Hurst parameter using a generalized method of moments. We establish their strong consistency under the two-step regime where the observation frequency $n$ and then the sample size $N$ tend to infinity, and prove their joint asymptotic normality when $H \in \big(\frac12, \frac34\big)$. Then, using a plug-in approach, we consistently estimate the random effects, and we study their asymptotic behavior under the same sequential asymptotic regime. Finally, we construct a nonparametric estimator for the distribution function of these random effects using a Lagrange interpolation at Chebyshev-Gauss nodes based method, and we analyze its asymptotic properties as both $n$ and $N$ increase. We illustrate the theoretical results through a numerical simulation framework. 
We further demonstrate the efficiency performance of the proposed estimators in an empirical application to crypto returns data, analyzing five major cryptocurrencies to uncover their distinct volatility structures and heterogeneous trend behaviors.\\


\section{Introduction}
In many practical contexts, different forms of randomness influence the dynamics of a system. Short-term fluctuations or noise with independent increments are commonly modeled using the Wiener process. However, this process does not capture the long-range dependence observed in a wide range of real-world systems. 
To address this limitation, fractional Brownian motion (fBm) has been proposed, as a generalization of standard Brownian motion that exhibits long-range dependence, making it well-suited for modeling dynamic systems characterized by memory effects and self-similarity. The statistical inference for stochastic differential equations (SDE) models driven by a standard Brownian motion has been extensively studied in many papers and summarized in several books, see e.g, \cite{bis08}, \cite{Kut04}, \cite{Pra99} and references cited therein. 
As for SDE models governed by an fBm, there has been considerable interest in studying problems of statistical inference for these models and how to estimate the unknown parameters on which the model depends, see e.g. \cite{Mish08}, \cite{Pr11} and references therein. 
Although models governed by the Wiener process or fBm capture either short-term noise or long-range dependence, respectively, many real-world systems exhibit both characteristics simultaneously. For instance, financial markets where asset prices exhibit both short-term volatility and long-term correlations. Similarly, models in telecommunications that have both rapid fluctuations and persistent network characteristics. Physiological signals often exhibit both fast, irregular fluctuations and underlying trends that persist over time.\\ 
In order to more accurately represent such phenomena and account for multiple noise sources, SDE models involving both a Wiener process and a fractional Brownian motion have been introduced. By combining these two processes within a single SDE framework, the resulting models leverage the strengths of each: the Wiener process captures short-term randomness, while fBm accounts for long-memory behavior. 
Several studies have investigated parameter estimation in stochastic models driven by mfBm. In particular, Prakasa Rao [\cite{Pra09}, \cite{Pra15a}, \cite{Pra15b}, \cite{Pra17}, \cite{Pra18}] developed an asymptotic theory for least squares and maximum likelihood estimators in various mfBm diffusion models. Mishura and Voronov \cite{Mish15} studied drift and diffusion estimation in mixed fractional Ornstein-Uhlenbeck processes. Mishra and Prakasa \cite{Mishr19} further extended these techniques to multidimensional and stochastic volatility settings; Mehrdoust et al. \cite{mehr20} considered a mixed-fractional Vasicek model and applied it to the pricing of Bermuda options on zero-coupon bonds, highlighting the relevance of mFBm in financial modeling.

In recent years, SDE models with random effects have been integrated into multiple research domains, including biology, pharmacokinetics, neuroscience, econometrics, and finance, as they enable simultaneous modeling of intrinsic system noise and subject-specific heterogeneity.  For instance, see \cite{delat12}, \cite{dit15}, \cite{don13} and \cite{pic10}. 
Nonparametric estimation for diffusion processes with random effects has been investigated in \cite{C13} and \cite{D14}. They developed kernel and deconvolution estimators of the common density of random effects and studied their asymptotic properties. 
In the case of fractional diffusion processes with random effects, El Omari et al. \cite{Om19} studied the properties of kernel and histogram estimators of the density of random effects and more recently, Chebli et al. \cite{chebli24} investigated a hybrid estimation approach to estimate the density of random effects by using a Bernstein polynomials based method. We also mention Prakasa Rao \cite{Pra20}, who studied the nonparametric estimation of the density of random effects in models governed by an SDE driven by a mixed fractional Brownian motion (in short mfBm).\\
Let $\left(\Omega,\mathcal{F},\left({\mathcal{F}}_t^i\right), \mathbb{P} \right)$ be a stochastic basis satisfying the usual conditions.\
We consider $N$ real valued stochastic processes $\left\{X^i_t,\ 0\leq t\leq T\right\}$, $i=1,\ldots, N$, with dynamics ruled by the following SDEs
\begin{equation}
	\begin{aligned}
		\label{B-S}
		\begin{cases} 
			dX^{i}_t &=\phi_{i}X^i_t dt +  X^i_t dM_t^{H,i},\\
			X^{i}_{0} &= x^i>0,
		\end{cases}
	\end{aligned}
\end{equation}
where $M^{H,i}=\left\{M^{H,i}_t,\ 0\leq t\leq T\right\}$, $i=1,\ldots,N$, are $N$ mutually independent mfBms such that $M^{H,i}_t:= \sigma B_t^i+\gamma B^{H,i}_t$, where $\sigma$ and $\gamma$ are positive constants to be estimated from discrete observations, $\left\{B^i_t,\ 0\leq t\leq T\right\}$ is a Brownian motion, and $\left\{B^H_t,\ 0\leq t\leq T\right\}$ is an independent fBm with unknown Hurst index $H\in (\frac{1}{2},1)$. The random effects 
$\phi_1,\ldots, \phi_N$ are $N$ unobserved independent and identically distributed (i.i.d.) random variables taking values in $\mathbb{R}$ and $x^i,\ i = 1, \ldots, N$, are known real values. Without loss of generality, we can assume that $x^1,\ldots,x^N=1$. We denote $F$ the common unknown distribution function of the random effects and $f$ their associated density function. The sequences $\{\phi_i,\ 1\le i\le N\}$ and $\{M^{H,i},\ 1\le i\le N\}$ are independent. 
Our main purpose is to estimate $F$, but first we should construct parametric estimators for the random effects $\phi_i$, $i=1,\ldots, N$, and for the fixed parameters $\sigma$, $\gamma$ and $H$ from the observations $\left\{X^i_t,\ 0\leq t\leq T\right\}$.
Since the Hurst parameter $H>\frac{1}{2}$, the sample trajectories of the fBm are almost surely Hölder continuous of any order less than $H$, and in particular, they are of bounded variation on compact intervals for $H>\frac{1}{2}$. Therefore, the stochastic integral $\displaystyle \int_0^t X_s^i d\,M^{H,i}_s$ can be interpreted pathwise in the Riemann–Stieltjes sense via Young integration, see \cite{nua}. Then, applying the pathwise fractional Itô formula, the model \eqref{B-S} admits a unique solution that is given for $i=1,\ldots,N$, by \begin{equation*}
	\begin{aligned} 
		X^i_t=\exp\left\{\left(\phi_i-\frac{1}{2}\sigma^2\right)t+ M^{H,i}_t\right\},\ 0\leq t\leq T.
	\end{aligned}
\end{equation*}
Denote for $ i=1,\ldots, N$ \begin{equation}\label{withtrend}
	\begin{aligned}
		Y^i_t:=\log\left(X^i_t\right) = \theta_i t + M_t^{H,i},\quad 0\leq t\leq T,
	\end{aligned}
\end{equation}
where $ \theta_i = \phi_i - \frac{\sigma^2}{2}.$ It's obvious that estimating parameters $\phi_i, \sigma^2, \gamma^2, H$ from model (\ref{B-S}) is equivalent to estimating them from model (\ref{withtrend}).\\
The process $\left\{Y^i_t,\ 0\leq t\leq T\right\}$ is called mixed fractional Brownian motion with trend.

Parameter estimation problem for this process in a classical framework without random effects, has attracted several researchers and works, including \cite{kuk22} where two types of estimators were proposed, one based on power variations and the other based on the ergodic theorem and the generalized method of moments which have allowed a simultaneous estimation for all the parameters. A similar simultaneous estimation was investigated in \cite{duf22}, where the maximum likelihood approach is combined with 
numerical optimization methods. A particular case of model (\ref{B-S}) where $\gamma=1$ and $H$, $\sigma$ are assumed to be known, was considered in \cite{cai16} where the authors constructed estimators for the drift parameter $\theta$. In \cite{Mish17}, drift parameter estimation was investigated for a more general model driven by a Gaussian process with stationary increments. 
More recently, \cite{Om24} considered a more general version of model (\ref{B-S}) without random effects, which he called $n$-th order mixed fractional Brownian motion with polynomial drift and constructed parametric estimators based on discrete observations for diffusion parameters and Hurst index, and on continuous observations for drift parameter, and studied asymptotic properties of all estimators.

To the best of our knowledge, the inference problem in stochastic diffusion models governed by mixed fractional Brownian motion and incorporating random effects has been addressed in only a few works. Prakasa Rao proposed in \cite{Pra20} a Kernel estimator for the density of the random effects based on continuous-time observations. In \cite{Pra21rand1}, Prakasa Rao studied the maximum likelihood estimation of the mean and variance of Gaussian random effects using discrete-time data. Finally, in  \cite{Pra21rand2}, the author investigated parametric estimation of the unknown parameters of a general random effects density. A key limitation in the previous studies is that the coefficients of the mfBm are not considered as unknown parameters to be estimated, but are constants equal to one. 
Our main contribution is to propose a complete hybrid estimation framework for the mixed fractional Black-Scholes model with random effects. We construct strongly consistent and asymptotically normal estimators for all fixed parameters $(H, \gamma^2, \sigma^2)$ and for the subject-specific random effects $\phi_i$. 
It is important to note that the theoretical derivation of the asymptotic normality for the parameter estimators requires the technical condition $H <\frac{3}{4}$, which is common in inference for fractional processes driven by mfBm. This condition ensures the finiteness of the asymptotic variance of the second-order moment statistics via the Breuer-Major theorem. The case $H \geq \frac{3}{4}$, which leads to non-Gaussian asymptotic limits (e.g., Rosenblatt distributions), is more complex and beyond the scope of this study. 
Furthermore, we introduce a novel nonparametric estimator for the distribution function $F$ of the random effects, based on Lagrange interpolation at Chebyshev-Gauss nodes. The motivation behind this choice is the fact that, in practical applications, for example, in censored or survival data in clinical or financial settings, the distribution of random effects often exhibits compact support. In such cases, classical kernel estimators admit considerable bias and variation near the boundaries of the compact support. In contrast, Lagrange interpolation at Chebyshev nodes respects the support's compactness and is shown to outperform classical kernel estimators, especially near the boundaries of the compact support. 

The paper is structured as follows. In Section \ref{sec2_chapi3}, we present some preliminaries on the properties of an mfBm and the procedure of nonparametric estimation based on Lagrange polynomials. Section \ref{sec3_chapi3} deals with the hybrid estimation procedure, which is split into two subsections. In the first, we construct parametric estimates of the random effects and the fixed parameters of our model in a discrete observations framework and the second is devoted to the nonparametric estimation of the distribution function $F$ and the asymptotic properties of the obtained estimator. Numerical simulations are presented in Section \ref{sec4_chapi3}. Section \ref{empiri} presents our empirical results, where we estimate the model on a pseudo-panel of five major cryptocurrencies using our developed estimation methods. Section \ref{conclusion_chapi3} contains some concluding remarks. To avoid interrupting the flow of the paper, the proofs of our results are relegated to Section \ref{sec6_chapi3}. We conclude this work with an Appendix where we recall some classical theorems used in this work.
\section{Preliminaries, notations and assumptions}
\label{sec2_chapi3}
This section is devoted to some notions that are mainly related to the mfBm and to the nonparametric functional estimation based on Lagrange polynomials.
\subsection{mfBm and related topics}\label{sub1}
Let $(\Omega,F, (\mathcal{F}_t)_{t\geq0}, P)$ be a stochastic basis satisfying the usual conditions. The natural filtration of a stochastic process is understood as the $\mathbb{P}$-completion of the filtration generated by this process. 
\begin{defi}\ \\
	A mixed fractional Brownian motion of parameters $a,\ b$, and $H$ is a process $M^H=\left\{M^H_t,\ t\geq 0\right\}$ defined by $M_t^H :=aB_t+bB_t^H,$ 
	where $\left\{B_t,\ t\geq0\right\}$ is a Brownian motion, and $\left\{B^H_t,\ t\geq0\right\}$ is an independent fractional Brownian motion with Hurst parameter $H\in (0,1)$ and $a$, $b$ are two real constants such that $\left(a,b\right)\neq \left(0,0\right)$.
\end{defi}
The mfBm $M^H$ satisfies the following properties:
\begin{itemize}
	\item $M^H$ is a centered Gaussian process and not a Markovian
	one.
	\item For all $t\in \mathbb{R}_{+}$, $\mathbb{E}\left(\left(M^H_t\right)^2\right)=a^2t+b^2t^{2H},$
	\item The covariance function of $M^H_t$ and $M^H_s$ for any $t,s\in \mathbb{R}_{+}$
	is given by
	$$Cov\left(M^H_t, M^H_s\right)=\frac{a^2}{2}\left(t+s-\lvert t-s\rvert\right)+\frac{b^2}{2}\left(t^{2H}+s^{2H}-\lvert t-s\rvert^{2H}\right).$$
	\item The increments of $M^H$ are positively correlated if $\frac{1}{2}< H <1$, uncorrelated if $H=\frac{1}{2}$ and negatively correlated if $0< H <\frac{1}{2}$,
	\item The increments of $M^H$ are long-range dependent if, and only if $H > \frac{1}{2}$.
	\item The increments of $M^H$ are stationary.
	\item  $M^H$ is a semimartingale in its own filtration if and only if either $H = \frac{1}{2}$ or $H\in (\frac{3}{4},1]$.
\end{itemize}
For further properties of the mfBm  and for details on the proofs of these properties, refer to \cite{ch01} and \cite{z06}.
\subsection{Approximation of a distribution function using Lagrange polynomials}
\label{sub2}
Let $Z_1,\ldots, Z_N$ be a sequence of i.i.d random variables with a common unknown distribution function $G$ supported on $[-1,1]$. 
The ordinary nonparametric estimator of order $m\geq 1$ for the distribution function $G$ obtained by using an approximation based on Lagrange polynomials with Chebyshev Gauss points, is defined for all $x\in [-1,1]$ as follows
\begin{equation}
	\begin{aligned}
		\label{lag}
		\tilde{G}_{m,N}(x)=\sum\limits_{j=1}^{m}\bar{G}_{N}\left(x_j\right)\mathcal{L}_{k}(x),
	\end{aligned}
\end{equation}
where \begin{itemize}
	\item $\bar{G}_N(y):=\dfrac{1}{N}\sum\limits_{i=1}^N \mathds{1}_{\left\{Z_i\leq y\right\}}$ is the empirical distribution function of the sequence $Z$.
	\item The points $x_j=\cos\left(\frac{(2k-1)\pi}{2m}\right)$, for all $k=1,\ldots, m$, denote the Chebyshev-Gauss nodes and they are the zeros of the Chebyshev polynomial \begin{eqnarray}T_m(x)=\cos\left(m \arccos\left(x\right)\right),\ x\in[-1,1].\label{Tm}\end{eqnarray} 
	\item $\mathcal{L}_j(x)=\prod\limits_{\underset{i\neq k}{i=1}}^m \frac{x-x_i}{x_j-x_i}$, $j=1,\ldots,m$ are the Lagrange polynomials.
\end{itemize}
For an in-depth reference on the properties of Lagrange polynomials with Chebyshev-Gauss points, we refer the reader to \cite{aus16}. The estimator (\ref{lag}) was proposed and studied in \cite{salima20} in a classical case without random effects.
In the case where the distribution function $G$ is supported on a compact interval $[a, b]$ with $a < b$, $Z$ can be transformed into a random variable $U$ supported on $[-1, 1]$, where $U = \frac{Z-\frac{(a + b)}{2}}{\frac{(b-a)}{2}}$. Transformations such as
$U = \dfrac{2Z}{1 +Z} - 1$ and $U = \left(2\pi \right)^{-1}\arctan\left(Z\right)$ can be used
to cover the cases of random variables with support $\mathbb{R}_{+}$
and $\mathbb{R}$ respectively. In the remainder of this paper, we assume that the distribution function $F$ of the random effects is supported on $\left[-1,1\right]$.\\
\subsection{Notation and assumptions}
Throughout this paper, we define $\log_2^{+}(x):=\mathds{1}_{\{x>0\}}\log_2 x$ and, for any bounded function $g$ on $[0,1]$, we set $\|g\|:=\sup_{x\in[-1,1]}|g(x)|$. The symbols $\overset{\mathrm{a.s.}}{\longrightarrow}$, $\overset{d}{\to}$, and $\xrightarrow{\mathbb{P}}$ denote almost sure convergence, convergence in distribution, and convergence in probability, respectively; finally, the superscript ${}^\top$ denotes the transpose of a vector.

The following assumptions are needed in the remainder of this paper.
\begin{assump}\ 
	\label{hypo}
	\begin{enumerate}
		\item[$(A1)$] For $i=1,\ldots,N$,   $\mathbb{E}[\phi_i^4] < \infty$.
		\item [$(A2)$] $F$ is of class $C^2$ on $[-1, 1]$.
		\item [$(A3)$] $f$ and $f'$ are bounded. 
	\end{enumerate} 
\end{assump}
\begin{com}\ 
	\begin{itemize}
		\item 
		Assumption (A1) is needed for studying the asymptotic properties of our model parameters via the delta method and the multivariate central limit theorem. It enables the control of squared terms that appear in second-order statistics, which we will use in our parametric estimation procedure. 

		\item Assumptions $(A2)$ and $(A3)$ are classical in nonparametric estimation theory using interpolation methods. The $C^2$ smoothness of $F$ allows for a second-order Taylor expansion, which is crucial for evaluating the bias of the Lagrange interpolant estimator $\widehat{F}_{m,n,N}$ in Proposition \ref{prop6_chap3}. The boundedness of $f$ and $f'$ is needed to control the behavior of the distribution function and the estimation error uniformly across its domain, including the boundaries.
	\end{itemize}
\end{com}
\section{Estimation procedure}
\label{sec3_chapi3}
In this section, we first investigate the parameter estimation of our model parameters and analyze their asymptotic properties. Then, in the second subsection, we construct a nonparametric estimator for the common distribution function of the random effects.
\subsection{Estimation of model parameters and random effects}
\par
This subsection is divided into two parts. First, we construct estimators for the global fixed parameter vector $\Theta = (H, \gamma^2, \sigma^2)$ using a generalized method of moments based on statistics averaged across all $N$ subjects. Second, we derive plug-in estimators for the subject-specific random effects $\phi_i$ based on the trajectory of each process.\\ 
Assume that the processes $\left\{Y^i_t,\ 0\leq t\leq T\right\}$ are observed at discrete times $t_k = kh,\ k = 0, \dots, n-1$, where $h > 0$ is a fixed time step and $n$ is the number of observations per subject. 
We define the discrete increments for each subject $i$ as follows
\begin{eqnarray}
	\Delta Y^i_k := Y^i_{(k+1)h} - Y^i_{kh} = \theta_i h + \sigma \Delta B^i_k + \gamma \Delta B^{H,i}_k= \theta_i h+\Delta M^{H,i}_k,\label{increm}
\end{eqnarray}
where 
\[\Delta M^{H,i}_k:=\sigma \Delta B^i_k + \gamma \Delta B^{H,i}_k,\quad
\Delta B^i_k = B^i_{(k+1)h} - B^i_{kh}, \quad \Delta B^{H,i}_k = B^{H,i}_{(k+1)h} - B^{H,i}_{kh}.
\]
The construction of all parameter estimators and the study of their asymptotic behavior are based on the following ergodic result, which will be proved in the proofs section.
\begin{lem}
	\label{lem_ergod}
	For each fixed $i=1,\ldots,N$, the process $\{\Delta Y^i_k-\theta_i h,\ 0\le k\le n-1\}$ is ergodic. 
\end{lem}
Based on this lemma, we present and study a natural estimator for the parameter $\hat{\theta}_i$ in the following proposition.
\begin{prop}\label{prop1_chap3}
	For each subject $i=1,\dots,N$, define
	\begin{eqnarray}
		\widehat{\theta}_i \;=\; \frac{1}{nh}\sum_{k=0}^{n-1}\Delta Y^i_k.\label{estheta}
	\end{eqnarray}
	Then, the following results hold.
	\begin{enumerate}
		\item[(i)] $ \widehat{\theta}_i$ is a strongly consistent estimator of $\theta_i$, that is,  
		$ \widehat{\theta}_i \;\xrightarrow{a.s.}\; \theta_i$, as $n\to\infty.$
		\item[(ii)] $ \widehat{\theta}_i$ is asymptotically normal and verifies
		\begin{eqnarray}
			n^{1-H}\big(\widehat{\theta}_i - \theta_i\big)\;\xrightarrow{d}\;\mathcal{N}\!\Big(0,\,\gamma^2 h^{2H-2}\Big).\label{asym_theta}
		\end{eqnarray}
	\end{enumerate}
\end{prop}

Now, based on the estimator $\hat{\theta}_i$ and the increments of $Y^i$, we introduce the following per-subject statistics that will be used in the construction of the parameter estimators.
\begin{align*}
	\xi^i_n := \dfrac{1}{n} \sum_{k=0}^{n-1} \left( \Delta Y^i_k  \right)^2, \quad
	\eta^i_n &:= \dfrac{1}{n} \sum_{k=0}^{n-1} \Delta Y^i_k \Delta Y^i_{k+1} , \quad
	\zeta^i_n := \dfrac{1}{n} \sum_{k=0}^{n-1} \left(Y^i_{(k+2)h}-Y^i_{kh}\right)\left(Y^i_{(k+4)h}-Y^i_{(k+2)h}\right).
\end{align*}
Since we are considering $N$ i.i.d subjects, we average the above statistics across all $N$ independent subjects as follows
\begin{align*}
	\bar{V}_{n,N} :=\frac{1}{N} \sum_{i=1}^N \widehat{\theta}_i^2,\quad
	\bar{\xi}_{n,N}:= \frac{1}{N} \sum_{i=1}^N \xi^i_n,\quad
	\bar{\eta}_{n,N} := \frac{1}{N} \sum_{i=1}^N \eta^i_n, \quad
	\bar{\zeta}_{n,N} := \frac{1}{N} \sum_{i=1}^N \zeta^i_n.
\end{align*}
The estimators of the global parameters and their asymptotic properties will be derived from the asymptotic behavior of the above empirical statistics. 
In the next proposition, we establish the almost sure convergence of the statistics 
\(
\bar{U}_{n,N} := \Big(\bar{V}_{n,N}, \bar{\xi}_{n,N}, \bar{\eta}_{n,N}, \bar{\zeta}_{n,N}\Big)^\top.
\)

\begin{prop}
	\label{cvg.as} 
	As $n,N \to \infty$ sequentially, the following almost sure convergence holds: 
	\[
	\bar{U}_{n,N} \xrightarrow{\text{a.s.}} U_{\infty},
	\]
	where $U_{\infty}:=\left(V_{\infty}, \xi_{\infty}, \eta_{\infty}, \zeta_{\infty}\right)^\top$ with
	\begin{eqnarray}
		V_{\infty} &:=& \mathbb{E}[\theta_i^2], \label{eq1}\\
		\xi_{\infty} &:=& \mathbb{E}[\theta_i^2] h^2 + \sigma^2 h + \gamma^2 h^{2H},\label{eq2}\\
		\eta_{\infty} &:=& \mathbb{E}[\theta_i^2] h^2 +  \gamma^2 h^{2H} \big(2^{2H - 1} - 1\big),\label{eq3}\\
		\zeta_{\infty} &:=& 4\mathbb{E}[\theta_{i}^2]h^2 + \gamma^2 h^{2H}\, 2^{2H}\big(2^{2H - 1} - 1\big). \label{eq4}
	\end{eqnarray}
\end{prop}
\begin{rem}
	The sequential asymptotic regime $n \to \infty$ first, then $N \to \infty$ is essential for the consistency and asymptotic normality of the subsequent global parameter estimators $\widehat{\Theta} = (\widehat{H}, \widehat{\gamma}^2, \widehat{\sigma}^2)^\top$. The reason is technical. Indeed, these estimators are constructed from the averaged statistics $(\bar{\xi}_{n,N},\bar{\eta}_{n,N},\bar{\zeta}_{n,N})$, whose limits require first letting $n \to \infty$ to use ergodicity within each subject’s trajectory and then $N \to \infty$ to apply the SLLN across subjects. If instead $n$ and $N$ increase simultaneously at comparable rates, the individual statistics $\xi^i_n$, $\eta^i_n$, and $\zeta^i_n$ may not be sufficiently close to their limits $\xi^i_\infty$, which introduces bias in the averaged statistics. Practically, the sequential regime corresponds to the realistic scenario where we collect dense observations ($n$ large) for each trajectory to accurately estimate its specific characteristics, and then aggregate information across many trajectories ($N$ large) to estimate population-level parameters.
	For this reason,  all asymptotic results are investigated under the sequential regime.
\end{rem}
The next proposition establishes the asymptotic normality of the empirical statistics.
\begin{prop}\label{prop3_chap3}
	Let $H\in\big(\frac12,\frac{3}{4}\big)$ and assume $(A1)$ hold. As $n,N \to \infty$ sequentially, we obtain
	\[
	\sqrt{N} \left( \bar{U}_{n,N} - U_\infty \right)
	\xrightarrow{d} \mathcal{N}(0,\, \Sigma),
	\]
	\begin{eqnarray}
		\Sigma:= \Var\big(\theta_i^2\big)\,\mathcal{H}\,\mathcal{H}^\top \quad \text{with}\quad \mathcal{H}=\big(1, h^2, h^2, 4h^2\big)^\top.\label{sigma1}\end{eqnarray}
	
\end{prop}

Using the empirical moments and their limits above,  in the next proposition, a strongly consistent estimator for the vector $\Theta$ is derived.
\begin{prop}
	\label{prop4_chap3}
	Under Assumption $(A1)$, $\widehat{\Theta}:=\left(\widehat{H},\widehat{\gamma}^2, \widehat{\sigma}^2\right)^\top$ is a strongly consistent estimator of $\Theta$, that is
	\[
	\widehat{\Theta}\xrightarrow{\text{a.s.}} \Theta
	\quad \text{as }\ n\to\infty,\ \mathrm{then}\ \ N\to \infty,
	\]
	where
	\begin{eqnarray}
		\widehat{H} &:=& \dfrac{1}{2} \log_2^+ \left( \dfrac{ \bar{\zeta}_{n,N} - 4h^2 \bar{V}_{n,N} }{ \bar{\eta}_{n,N} - h^2 \bar{V}_{n,N} } \right),\label{hhat}
		\\
		\widehat{\gamma}^2 &:=& \dfrac{ \bar{\eta}_{n,N} - h^2 \bar{V}_{n,N} }{ h^{2\widehat{H}} \left(2^{2\widehat{H} - 1} - 1\right) },\label{gam} \\
		\widehat{\sigma}^2 &:=& \dfrac{ \bar{\xi}_{n,N}- h^2 \bar{V}_{n,N} - \widehat{\gamma}^2 h^{2\widehat{H}} }{ h }\label{sig}.
	\end{eqnarray}
\end{prop}
The next theorem concerns the asymptotic normality of $\widehat{\Theta}_N$.
		\begin{theor}\label{theo1_chap3}
			Let assumption $(A1)$ hold. Then, for $H\in\big(\tfrac12,\tfrac34\big)$ and as $n,N \to \infty$ sequentially, the estimator $\widehat{\Theta}_N$ is asymptotically normal, that is
			\[
			\sqrt{N}\,\big(\widehat{\Theta}_N-\Theta\big)\ \xrightarrow{d}\ \mathcal{N}\!\left(0,\;\mathcal{C}\right),
			\]
			where the asymptotic covariance matrix $\mathcal{C}$ is given by $\mathcal{C}:=J\,\Sigma\,J^{\!\top}$ with $\Sigma$ is the matrix \eqref{sigma1} and the Jacobian matrix $J$ is defined as follows
			\[
			J:=\left.\frac{\partial \widehat{\Theta}}{\partial \bar U_{n,N}}\right|_{U_\infty}
			=\begin{pmatrix}
				\left.\dfrac{\partial \widehat H}{\partial \bar V_{n,N}}\right|_{U_\infty} &
				\left.\dfrac{\partial \widehat H}{\partial \bar \xi_{n,N}}\right|_{U_\infty} &
				\left.\dfrac{\partial \widehat H}{\partial \bar \eta_{n,N}}\right|_{U_\infty} &
				\left.\dfrac{\partial \widehat H}{\partial \bar \zeta_{n,N}}\right|_{U_\infty} \\[2mm]
				\left.\dfrac{\partial \widehat\gamma^{2}}{\partial \bar V_{n,N}}\right|_{U_\infty} &
				\left.\dfrac{\partial \widehat\gamma^{2}}{\partial \bar \xi_{n,N}}\right|_{U_\infty} &
				\left.\dfrac{\partial \widehat\gamma^{2}}{\partial \bar \eta_{n,N}}\right|_{U_\infty} &
				\left.\dfrac{\partial \widehat\gamma^{2}}{\partial \bar \zeta_{n,N}}\right|_{U_\infty} \\[2mm]
				\left.\dfrac{\partial \widehat\sigma^{2}}{\partial \bar V_{n,N}}\right|_{U_\infty} &
				\left.\dfrac{\partial \widehat\sigma^{2}}{\partial \bar \xi_{n,N}}\right|_{U_\infty} &
				\left.\dfrac{\partial \widehat\sigma^{2}}{\partial \bar \eta_{n,N}}\right|_{U_\infty} &
				\left.\dfrac{\partial \widehat\sigma^{2}}{\partial \bar \zeta_{n,N}}\right|_{U_\infty}
			\end{pmatrix}\!,
			\]
			where the derivatives of $\widehat{H}$ are given as follows 
			\[
			\frac{\partial \widehat H}{\partial V}\Big|_{U_\infty}
			=\frac{h^2}{2\ln 2}\,\frac{2^{2H}-4}{2^{2H}\gamma^2 D_H},\qquad
			\frac{\partial \widehat H}{\partial \xi}\Big|_{U_\infty}=0,\qquad
			\frac{\partial \widehat H}{\partial \eta}\Big|_{U_\infty}=-\frac{1}{2\ln 2}\,\frac{1}{B},\qquad
			\frac{\partial \widehat H}{\partial \zeta}\Big|_{U_\infty}=\frac{1}{2\ln 2}\,\frac{1}{2^{2H} B},
			\]
			the derivatives of $\widehat{\gamma}^2$ are expressed as 
			\[
			\begin{aligned}
				\frac{\partial \widehat\gamma^{2}}{\partial V}\Big|_{U_\infty}
				&= -\,\frac{h^2}{D_H} + \frac{h^2}{2\ln 2}\!\left(\frac{4}{2^{2H}}-1\right)\frac{D'_H}{D_H^{2}},\qquad
				\frac{\partial \widehat\gamma^{2}}{\partial \xi}\Big|_{U_\infty}=0,\\[1mm]
				\frac{\partial \widehat\gamma^{2}}{\partial \eta}\Big|_{U_\infty}
				&= \frac{1}{D_H} + \frac{D'_H}{2\ln 2}\,\frac{1}{D_H^{2}},\qquad
				\frac{\partial \widehat\gamma^{2}}{\partial \zeta}\Big|_{U_\infty}
				= -\,\frac{D'_H}{2\ln 2}\,\frac{1}{2^{2H} D_H^{2}},
			\end{aligned}
			\]
			and finally the derivatives of $\widehat{\sigma}^2$ are given by
			\[
			\begin{aligned}
				\frac{\partial \widehat\sigma^{2}}{\partial V}\Big|_{U_\infty}
				&= -\,h\;+\;h^{\,2H-1}\!\Bigg[
				\frac{h^{2}}{D_H}
				-\frac{h^{2}}{2\ln 2}\!\left(1-\frac{4}{2^{2H}}\right)
				\left(\!-\frac{D'_H}{D_H^{2}}+\frac{2\ln(h)}{D_H}\!\right)
				\Bigg],\\[1mm]
				\frac{\partial \widehat\sigma^{2}}{\partial \xi}\Big|_{U_\infty}&=\frac{1}{h},\\[1mm]
				\frac{\partial \widehat\sigma^{2}}{\partial \eta}\Big|_{U_\infty}
				&= -\,h^{\,2H-1}\,\frac{1}{D_H}
				\left[\,1-\frac{2\ln(h)}{2\ln 2}+\frac{D'_H}{2\ln 2}\,\frac{1}{D_H}\right],\\[1mm]
				\frac{\partial \widehat\sigma^{2}}{\partial \zeta}\Big|_{U_\infty}
				&= \;h^{\,2H-1}\,2^{-2H}\!
				\left[\,\frac{D'_H}{2\ln 2}\,\frac{1}{D_H^{2}}
				-\frac{2\ln(h)}{2\ln 2}\,\frac{1}{D_H}\right],
			\end{aligned}
			\]
			with 
			\(
			B:=\eta_\infty-h^2V_\infty=\gamma^2 h^{2H}\!\big(2^{2H-1}-1\big),\quad
			A:=\zeta_\infty-4h^2V_\infty=2^{2H} B,\quad
			D_H:=h^{2H}\!\big(2^{2H-1}-1\big),\) \\ and \(
			D'_H:=h^{2H}\!\Big[2\ln(h)\!\big(2^{2H-1}-1\big)+ (2\ln 2)\,2^{2H-1}\Big].
			\)
		\end{theor}
		\begin{rem}\label{Hrange}
			The restriction $H \in \left(\frac{1}{2}, \frac{3}{4}\right)$ ensures that the statistics vector $U_{n,N}$ satisfies a Gaussian central limit theorem under the scaling $\sqrt{N}$. This stems from the Breuer-Major theorem (Theorem  \ref{bru_m}), which requires the condition
			\begin{equation}
				\sum_{k \in \mathbb{Z}} |\rho(k)|^2 < \infty,
			\end{equation}
			where $\rho(k)$ is the autocovariance function defined in \eqref{ro}. Using the asymptotic behavior of $\rho_H(k)$, which is given by $\rho_H(k) \sim h^{2H}H(2H-1)k^{2H-2}$ as $k \to \infty$, we obtain
			\begin{equation}
				\sum_{k=1}^\infty |\rho_H(k)|^2 \sim C h^{4H} \sum_{k=1}^\infty k^{4H-4},
			\end{equation}
			which converges if and only if $4H - 4 < -1$, i.e., $H < \frac{3}{4}$. At the critical value $H = \frac{3}{4}$, the asymptotic behavior of functionals of fractional Gaussian noise typically changes: for $H < \frac{3}{4}$, one obtains a Gaussian limit under standard normalizations, whereas for $H \geq \frac{3}{4}$ the limit becomes non-Gaussian and typically involves Rosenblatt-type distributions (see for instance \cite{dob}, \cite{Mish18}, \cite{rosen}, \cite{taqqu}). The case $H \geq \frac{3}{4}$ requires different techniques involving multiple Wiener-Itô integrals and Malliavin calculus, which is beyond the scope of this paper.
		\end{rem}
		Now, relying on the previous construction and study of the plug-in estimates $\widehat{\Theta}=\left(\widehat{H}, \widehat{\gamma}^2,\widehat{\sigma}^2\right)^\top$, we can recover the unobserved subject-specific random effects $\phi_i:=\theta_i + \frac{1}{2} \sigma^2$, for $i=1,\ldots,N$.\\
				Substituting $\theta_i$ and $\sigma^2$ by their estimates $\widehat{\theta}_i$ and $\widehat{\sigma}^2$ investigated in \eqref{estheta} and \eqref{sig} respectively, a plug-in estimator of the random effect $\phi_i$ is given for each $i=1,\ldots,N$ by \begin{eqnarray}
					\widehat{\phi}_i := \widehat{\theta}_i + \tfrac{1}{2} \widehat{\sigma}^2.\label{reffhat}
				\end{eqnarray}
				
				\begin{rem}\label{dependence}
					The estimators $\widehat\phi_1, \ldots, \widehat\phi_N$ are not independent because they all use the same plug-in estimator
					$\widehat\sigma^2$. Under our sequential regime, 
					we show that $\max\limits_{1\le i\le N}|\widehat\phi_i-\phi_i|\overset{a.s}{\longrightarrow}0$ (see equation \eqref{maxphi}).\\ Hence, the dependence does not affect the uniform consistency. For distributional limits, the common error $\frac{1}{2}(\widehat\sigma^2-\sigma^2)=O_{\mathbb{P}}(N^{-\frac12})$ must be
					controlled at the corresponding rate; this is done explicitly in the proofs of Theorems \ref{theo3_chap3}-\ref{theo4_chap3} and Proposition \ref{prop7_chap3}).
				\end{rem}
				%
				%
				The following Proposition provides the strong consistency and asymptotic normality of $\widehat{\phi}_i, i = 1, \dots, N$, for the sequential regime ($n \to \infty$ then $N \to \infty$).
				\begin{prop}\ \label{prop5_chap3}
					Let $H\in\Big(\frac12, \frac34\Big)$. Under assumption $(A1)$, the following results hold.
					\begin{enumerate}
						\item As $n,N \to \infty$ sequentially, the estimator $\widehat{\phi}_i$ is strongly consistent, that is, for each fixed $i=1,\ldots,N$, 
						$\widehat{\phi}_i\ \xrightarrow{a.s} \ \phi_i.$
							\item As $n,N \to \infty$ sequentially, we obtain
							\[
							\sqrt N (\widehat\phi_i-\phi_i) \xrightarrow{d} \mathcal N\Big(0,\frac14 v_\sigma^2\Big),
							\]
							where $v_\sigma^2$ is the asymptotic variance of $\widehat{\sigma^2}$.
						\end{enumerate}
					\end{prop}

					\vspace{1em}
					
					In the next subsection, we proceed to estimate the distribution \( F \) of random effects based on the estimated random effects \( \widehat{\phi}_{1},\ldots, \widehat{\phi}_{N}\).
					\subsection{Nonparametric estimation of random effects distribution function}
					For the random effects distribution function $F$, the estimator $(\ref{lag})$ has the following form:
					\begin{equation*}
						\begin{aligned}
							\label{Fhat}
							\tilde{F}_{m,N}(x)=\sum\limits_{j=1}^{m}\bar{F}_{N}\left(x_j\right)\mathcal{L}_{j}(x),\ \forall x\in [-1,1],
						\end{aligned}
					\end{equation*}
					where $\bar{F}_N(y)=\frac{1}{N}\sum\limits_{i=1}^N \mathds{1}_{\left\{\phi_i\leq y\right\}}$ denotes the empirical distribution function of the random effects $\phi_i$.
					As we can observe, this estimator depends on $\phi_i$ which are not observed. Hence, it is natural to replace them with their estimates to obtain a suitable calculable estimator of $F$. 
					Substituting the random effects $\phi_i$ by their estimators investigated in (\ref{reffhat}), we obtain the following estimator of $F$,
					\begin{equation}\begin{aligned}\label{F_phihat}\widehat{F}_{m,n,N}(x)=\sum\limits_{j=1}^{m}\widehat{F}_{N}\left(x_j\right)\mathcal{L}_{j}(x),\ \forall x\in [-1,1],\end{aligned}\end{equation}
					where $\widehat{F}_N(y)=\frac{1}{N}\sum\limits_{i=1}^N \mathds{1}_{\left\{\widehat{\phi}_{i}\leq y\right\}}$ is the empirical distribution function of the estimators $\widehat{\phi}_{i}$. The remainder of this subsection is devoted to studying the asymptotic properties of the estimator $\widehat{F}_{m,n,N}$. \\
					The following proposition sets forward the asymptotic bias, variance and mean squared error (in short MSE) of $\widehat{F}_{m,n,N}$.
					\begin{prop}
						\label{prop6_chap3}
						Let assumptions $(A1)-(A3)$ hold. Then, for $x\in[-1,1]$,
						\begin{enumerate} \item the asymptotic bias of $\widehat{F}_{m,n,N}$ is given by
							\begin{equation}
								\label{bi}
								\lim\limits_{n,N\to \infty}\mathrm{Bias}\left(\widehat{F}_{m,n,N}(x)\right)=\frac{\pi}{m^2}T_m(x) \mu(x)+O(m^{-2}),
							\end{equation}
							where $T_m$ is defined \eqref{Tm} and $\mu(x):=\frac{1}{4}(x-1)f'(x)-\frac{1}{2}f(x)$.
							\item the asymptotic variance of $\widehat{F}_{m,n,N}$ is given by
							\begin{equation}
								\label{vari}
								\lim\limits_{n,N\to \infty} Var\left(\widehat{F}_{m,n,N}(x)\right)=N^{-1} \sigma_F^{2}(x)+O\left(N^{-1}m^{-\frac 12}\right),
							\end{equation}
							where $\sigma_F^2(x)=F(x)\left(1-F(x)\right)$.
							\item the asymptotic mean squared error of $\widehat{F}_{m,n,N}$, denoted by $AMSE\left(\widehat{F}_{m,n,N}(x)\right)$, is given by 
							\begin{eqnarray} \mathrm{AMSE}\left(\widehat{F}_{m,n,N}(x)\right)=\dfrac{\pi}{m^4}C_1(x) + N^{-1} \sigma_F^{2}(x)+ O(m^{-4}) + O\left(N^{-1}m^{-\frac 12}\right), \label{limse} \end{eqnarray}
							where $C_1(x):=T^2_m(x) \mu^2(x)$.
						\end{enumerate}	
					\end{prop}
					
					
					
					The next proposition deals with the uniform strong consistency of $\widehat{F}_{m,n,N}$.
					\begin{theor}
						\label{theo2_chap3}
						Let assumptions $(A1)-(A3)$ hold and $m=m_N$ satisfy $m_N\to\infty$ and $\displaystyle\frac{\log m_N}{m_N^{2}}\to 0$ as $N\to\infty$. Then, 
						Under the sequential regime $n\to\infty$ first, then $N\to\infty$, there exists a sequence $n=n_N\to\infty$ such that
						\[
						\bigl\|\widehat F_{m,n,N}-F\bigr\| \overset{a.s}{\longrightarrow} 0,\qquad N\to\infty.
						\]
					\end{theor}
					Now we turn our attention to the asymptotic normality of the estimator.
					\begin{theor}
						\label{theo3_chap3}
						Assume $(A1)-(A3)$. Let $m=m_N$ satisfy
						\(
						N^{-\frac 12}\,\frac{\log m}{m^2}\;\longrightarrow\;0
						\quad\text{as }N\to\infty .
						\)
						Then, under the sequential asymptotic regime $n\to\infty$ followed by $N\to\infty$, for every fixed
						$x\in[-1,1]$,
						\[
						\sqrt N\bigl(\widehat F_{m,n,N}(x)-F(x)\bigr)
						\;\xrightarrow{d}\;
						\mathcal N\!\left(0,\sigma_F^2(x)\right),
						\qquad
						\sigma_F^2(x)=F(x)\bigl(1-F(x)\bigr).
						\]
					\end{theor}
					In the following theorem, we establish an explicit non-asymptotic upper bound for the mean squared error of $\hat{F}_{m,n,N}$, quantifying the contribution of the interpolation error, the finite-sample effect, and the error induced by estimating the random effects.
					\begin{theor}\label{theo4_chap3}
						Assume (A1)-(A3) are satisfied and let $m = m_N$ and $n = n_N$. 
						Let $\varepsilon_N \to 0$ be a sequence satisfying $$\max\limits_{1 \leq i \leq N} |\widehat{\phi}_i - \phi_i| \leq \varepsilon_N,\quad \text{a.s.}$$ 
						Then
						\begin{equation}\label{bound_F}
							\mathbb{E}\left[\|\widehat{F}_{m,n,N} - F\|^2\right] \leq \frac{C_1 (\log m_N)^2}{m_N^4} + \frac{C_2}{N} + C_3 \varepsilon_N^2 + o\left(\frac{1}{m_N^4} + \frac{\varepsilon_N}{N}\right),
						\end{equation}
						where $C_1$, $C_2$, $C_3$ are positive constants depending on $\|f\|$, $\|f'\|$ and $h$.
					\end{theor}
					\begin{rem}
						A natural way to balance the contributions of the non-asymptotic bound terms is to choose 
						$m_N\sim N^{\frac 14}$, yielding
						$\frac{(\log m_N)^2}{m_N^4}\asymp\frac{(\log N)^2}{N}$. Moreover, under the sequential regime $n\to\infty$ followed by $N\to\infty$, Proposition \ref{prop5_chap3} and Theorem \ref{theo1_chap3} ensure that
						\[
						\hat\phi_i-\phi_i
						=O_{\mathbb{P}}\big(n_N^{H-1}\big)+O_{\mathbb{P}}\big(N^{-\frac12}\big).
						\]
						Then, by choosing $n_N=N^2$ and under the condition $H<\frac34$, we obtain 
						$\varepsilon_N=O_{\mathbb{P}}\big(N^{-\frac12}\big)$, so that
						$C_3\varepsilon_N^2\asymp N^{-1}$. 
						With these choices, 
						the three leading terms are all of order $N^{-1}$. Consequently, there exists a constant $C>0$ such that, for $N$ large enough,
						\[
						\mathbb{E}\big\|\hat F_{m_N,n_N,N}-F\big\|^2
						\;\le\;
						C\,\frac{(\log N)^2}{N}.
						\]
						Hence, our non-asymptotic bound demonstrates that the estimator $\hat F_{m_N,n_N,N}$, which is based on plug-in estimators of the random effects, achieves, up to a logarithmic factor, the same convergence rate as the empirical c.d.f $\bar F_N$ would if the random effects were directly observed.  
					\end{rem}

					In the next proposition, we show that, under certain conditions, the proposed Lagrange estimator satisfies the Chung-Smirnov property, which characterizes its extreme fluctuations about $F$ as the number of subjects $N$ increases. This property was investigated for the empirical c.d.f by \cite{chu} and \cite{smir}.
					\begin{prop}
						Let assumptions $(A1)$-$(A3)$ hold. Assume that $n$ and $N$ tend to infinity sequentially and let $m=m_N\to\infty$ be such that
						$\left(\frac{2N}{\log\log N}\right)^{\frac 12}\frac{\log m}{m^2}\longrightarrow 0.
						$
						Then $\widehat{F}_{m,n,N}$ satisfies the Chung-Smirnov property, that is 
						\begin{eqnarray*}
							\underset{N\to \infty}{\lim \sup} \left(\dfrac{2N}{\log\log N}\right)^{\frac{1}{2}} \underset{x\in[-1,1]}{\sup}\left|\widehat{F}_{m,n,N}(x)-F(x)\right| \leq 1, \ \text{a.s}.
						\end{eqnarray*}
						\label{prop7_chap3}
					\end{prop}

					In order to illustrate the theoretical results discussed above, we present in the following section some numerical simulations and graphical results implemented using Python.
					\section{Numerical simulations}\label{sec4_chapi3}
					In this section, we evaluate the performance of our hybrid estimation method based on simulated data. In order to test our parametric estimation method, we use a Monte-Carlo simulation with 50 iterations. The increments $\Delta Y^i$, $i=1,\ldots, N$ are generated using a combination of the Euler-Maruyama scheme for the Brownian motion increments and Cholesky decomposition for the fBm increments. The true parameters are fixed as $H=0.7$, $\gamma=0.5$ and $\sigma=0.2$, and the random effects $\phi_i$ are drawn from each example of the following examples of distributions
					\begin{itemize}
						\item[(1)] Beta distribution $\mathcal{B}(2,2)$,
						\item[(2)] Gamma distribution $\mathcal{G}(2,1)$,
						\item[(3)] Gaussian distribution $\mathcal{N}(0.5,0.25)$,
						\item[(4)] Gaussian mixture distribution $0.5*\mathcal{N}(-2,1) + 0.5* \mathcal{N}(3,0.5)$.
					\end{itemize}
					Estimation methods are tested for $N\in\{100, 250, 500\}$ subjects and $n\in\{250, 500, 1000\}$ observations per subject. 
					The simulated results involve the mean and standard deviation (S.dev) of the estimators $\widehat{H},\ \widehat{\gamma}^2$ and $\widehat{\sigma}^2$, as well as the random effect estimator $\widehat{\phi}_i$ and are presented in Table \ref{tab2}.
					
					\begin{table}[ht]
						\centering
						\caption{Parametric estimation results with true values $H = 0.7$, $\gamma^2 = 0.25$ and $\sigma^2 = 0.04$.}
						\resizebox{\textwidth}{!}{%
							\begin{tabular}{ccccccccccc}
								\hline
								$F$ & $(N,n)$ 
								& \multicolumn{2}{c}{$\widehat{H}$} 
								& \multicolumn{2}{c}{$\widehat{\gamma}^2$} 
								& \multicolumn{2}{c}{$\widehat{\sigma}^2$} 
								& \multicolumn{3}{c}{$\widehat{\phi}_{i}$} \\
								\cline{3-11}
								& & Mean & S.dev & Mean & S.dev & Mean & S.dev & True Mean & Mean($\widehat{\phi}_{i}$) & S.dev \\
								\hline
								\textbf{True Value} & & \textbf{0.7} & - & \textbf{0.25} & - & \textbf{0.04} & - & - & - & - \\
								\hline
								\multirow{3}{*}{(1)}  
								& (100, 250)  & 0.7009 & 0.0093 & 0.2159 & 0.0150 & 0.0442 & 0.0021 & 0.5077 & 0.5085 & 0.0426 \\
								& (250, 500) & 0.6997 & 0.0094 & 0.2290 & 0.0095 & 0.0422 & 0.0019 & 0.5066 & 0.5073 & 0.0175 \\
								& (500, 1000) & 0.6993 & 0.0106 & 0.2361 & 0.0092 & 0.0412 & 0.0023 & 0.5080 & 0.5078 & 0.0092 \\
								\hline
								\multirow{3}{*}{(2)} 
								& (100, 250)  & 0.6985 & 0.0092 & 0.2155 & 0.0119 & 0.0435 & 0.0024 & 1.9192 & 1.9200 & 0.0409 \\
								& (250, 500) & 0.6955 & 0.0098 & 0.2297 & 0.0100 & 0.0422 & 0.0025 & 2.0624 & 2.0641 & 0.0216 \\
								& (500, 1000) & 0.7013 & 0.0104 & 0.2382 & 0.0023 & 0.0401 & 0.0003 & 1.9834 & 1.9818 & 0.0123 \\
								\hline
								\multirow{3}{*}{(3)} 
								& (100, 250)  & 0.7003 & 0.0110 & 0.2170 & 0.0131 & 0.0438 & 0.0028 & 0.4681 & 0.4630 & 0.0409 \\
								& (250, 500) & 0.6990 & 0.0098 & 0.2282 & 0.0078 & 0.0420 & 0.0025 & 0.5188 & 0.5165 & 0.0187 \\
								& (500, 1000) & 0.6988 & 0.0104 & 0.2359 & 0.0094 & 0.0411 & 0.0023 & 0.5234 & 0.5255 & 0.0100 \\
								\hline
								\multirow{3}{*}{(4)} 
								& (100, 250)  & 0.7003 & 0.0095 & 0.2169 & 0.0130 & 0.0438 & 0.0028 & 0.3935 & 0.3884 & 0.0409 \\
								& (250, 500) & 0.6990 & 0.0098 & 0.2283 & 0.0079 & 0.0419 & 0.0025 & 0.4877 & 0.4854 & 0.0187 \\
								& (500, 1000) & 0.6988 & 0.0104 & 0.2359 & 0.0094 & 0.0411 & 0.0023 & 0.5045 & 0.5025 & 0.0100 \\
								\hline
							\end{tabular}
						}
						\label{tab2}
					\end{table}
					Using the simulated estimators, we numerically compare the proposed Lagrange estimator \eqref{F_phihat} with the standard kernel estimator of the random effects distribution function, which is defined by \begin{eqnarray}\label{kern_phi}\widehat{F}_h(x):=\dfrac{1}{N}\sum\limits_{i=1}^N K\left(\dfrac{x-\widehat{\phi}_i}{h}\right), \quad \forall \ x\in\mathbb{R},\end{eqnarray} where $K$ is given by $K(y):=\displaystyle\int_{-\infty}^{y} k(u)\,du$ and $k$ is a kernel function satisfying $$\displaystyle\int_{-\infty}^{+\infty}k(u)\,du=1,\ \displaystyle\int_{-\infty}^{+\infty}u\,k(u)\,du=0\quad \text{and} \ \displaystyle\int_{-\infty}^{+\infty}u^2\,k(u)\,du\neq 0.$$ 
					Here, the Gaussian kernel is chosen. For each distribution function example from (1) to (4), we use the relevant transformation introduced in Section \ref{sub2} and test both estimators under two sample size regimes: $N=100$ subjects with $n=250$ observations per subject and $N=250$ subjects with $n=500$ observations. For the choice of the interpolation order$m$, we use $K$-fold cross-validation. 
					We first partition the estimated random effects $\{\widehat{\phi}_i, i = 1,\ldots,N\}$ into $K$ folds of approximately equal size, with index sets $I_1,\ldots,I_K$. For each candidate $m \in \mathcal{M}$, and for each fold $k$:
					\begin{itemize}
						\item we construct the Lagrange estimator $\widehat{F}^{(-k)}_{m}$ using all data points with indices in $\{1,\ldots,N\}\setminus I_k$.
						\item we evaluate its prediction error on the validation fold via
						\[
						\mathrm{VE}_k(m)
						=
						\frac{1}{|I_k|}
						\sum_{i\in I_k}
						\bigl[\widehat{F}^{(-k)}_{m}(\widehat{\phi}_i) - \widehat{F}^{(k)}(\widehat{\phi}_i)\bigr]^2,
						\]
						where $\widehat{F}^{(k)}$ denotes the empirical c.d.f.\ based only on observations in $I_k$.
					\end{itemize}
					The cross-validation score for $m$ is then
					\[
					\mathrm{VE}(m) = K^{-1}\sum_{k=1}^K \mathrm{VE}_k(m),
					\]
					and we select
					\[
					m_{\mathrm{CV}} = \arg\min_{m\in\mathcal{M}} \mathrm{VE}(m).
					\]
					In the presented simulations, we take $K=5$ and $\mathcal{M}=\{5,\ldots,20\}$, 
					and we denote
					the selected value by $m_{\text{opt}}$. The final Lagrange estimator $\widehat{F}_{m,n,N}$
					is then constructed using $m = m_{\text{opt}}$.
					
					The plots of both estimators are presented in the figures that follow. 
					\begin{figure}[ht]
						\centering
						\begin{subfigure}[b]{0.48\textwidth}
							\includegraphics[width=\linewidth]{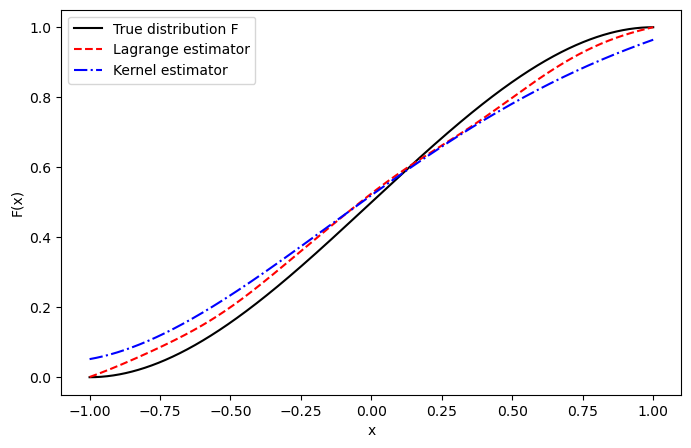}
							\caption{$N=100$, $n=250$}
						\end{subfigure}
						\hfill 
						\begin{subfigure}[b]{0.48\textwidth}
							\includegraphics[width=\linewidth]{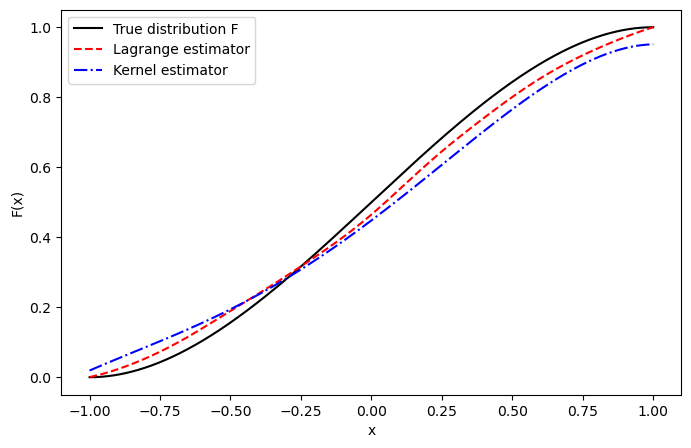}
							\caption{$N=250$, $n=500$}
						\end{subfigure}
						\caption{Qualitative comparison between the proposed estimator $\widehat{F}_{m,n,N}$ and the Kernel estimator $\widehat{F}_h$ for the distribution $\mathcal{B}(2,2)$. }
						\label{fig1}
					\end{figure}
					\begin{figure}[ht]
						\centering
						\begin{subfigure}[b]{0.48\textwidth}
							\includegraphics[width=\linewidth]{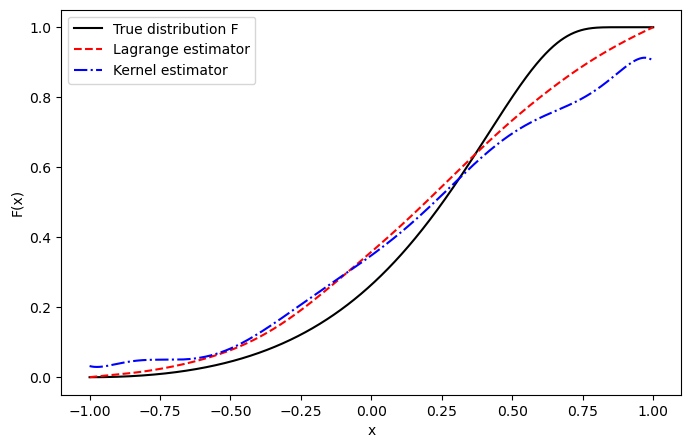}
							\caption{$N=100$, $n=250$}
						\end{subfigure}
						\hfill 
						\begin{subfigure}[b]{0.48\textwidth}
							\includegraphics[width=\linewidth]{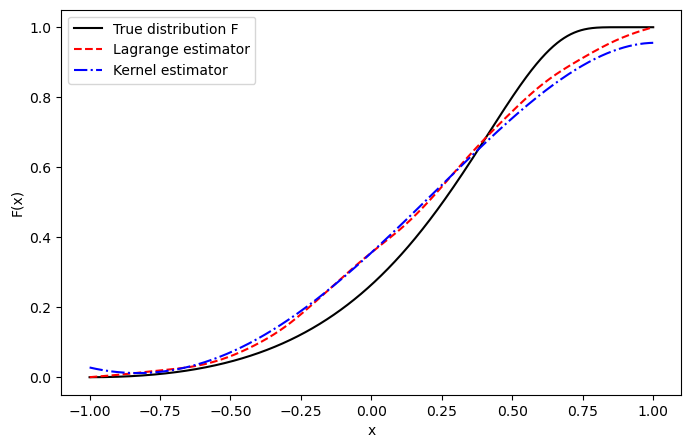}
							\caption{ $N=250$, $n=500$}
						\end{subfigure}
						\caption{Qualitative comparison between the proposed estimator $\widehat{F}_{m,n,N}$ and the Kernel estimator $\widehat{F}_h$ for the distribution $\mathcal{G}(2,1)$. }
						\label{fig2}
					\end{figure}
					\begin{figure}[ht]
						\centering
						\begin{subfigure}[b]{0.48\textwidth}
							\includegraphics[width=\linewidth]{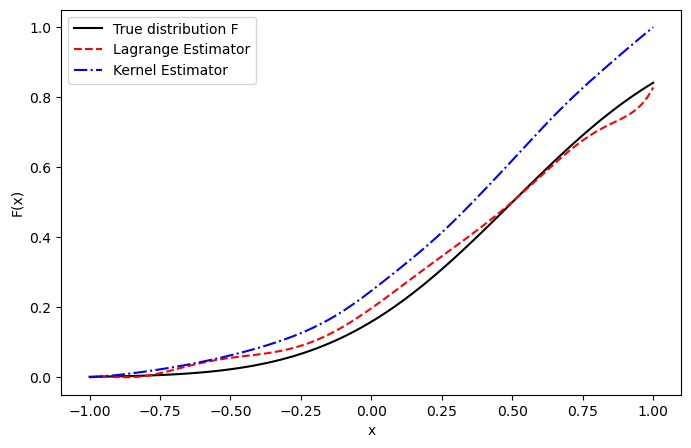}
							\caption{ $N=100$, $n=250$}
						\end{subfigure}
						\begin{subfigure}[b]{0.48\textwidth}
							\includegraphics[width=\linewidth]{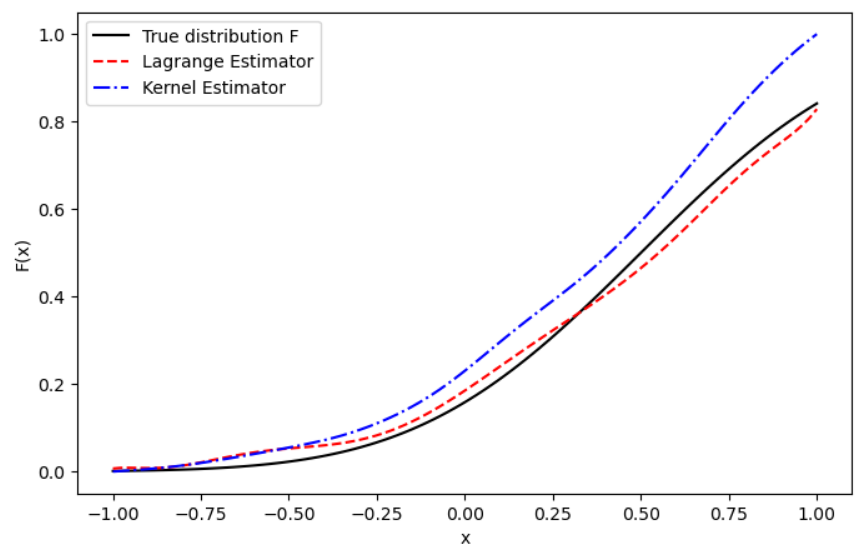}
							\caption{$N=250$, $n=500$}
						\end{subfigure}
						\caption{Qualitative comparison between the proposed estimator $\widehat{F}_{m,n,N}$ and the Kernel estimator $\widehat{F}_h$ for the distribution $\mathcal{N}(0.5,0.25)$. }
						\label{fig3}
					\end{figure}
					\begin{figure}[ht]
						\centering
						\begin{subfigure}[b]{0.48\textwidth}
							\includegraphics[width=\linewidth]{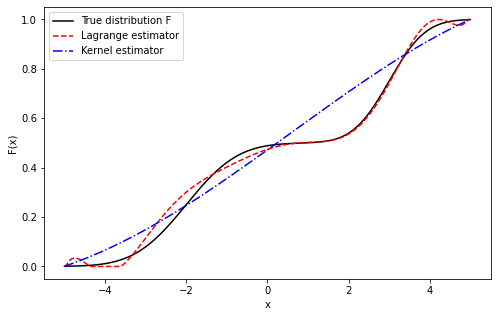}
							\caption{ $N=100$, $n=250$}
						\end{subfigure}
						\hfill 
						\begin{subfigure}[b]{0.48\textwidth}
							\includegraphics[width=\linewidth]{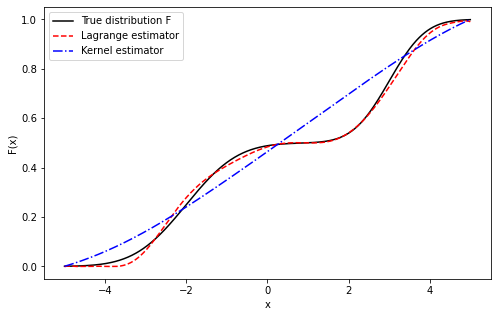}
							\caption{$N=250$, $n=500$}
						\end{subfigure}
						\caption{Qualitative comparison between the proposed estimator $\widehat{F}_{m,n,N}$ and the Kernel estimator $\widehat{F}_h$ for the distribution $0.5*\mathcal{N}(-2,1) + 0.5* \mathcal{N}(3,0.5)$.}
						\label{fig4}
					\end{figure}
					
					In order to numerically compare the performance of $\widehat{F}_{m,n,N}$ and $\widehat{F}_h$, we compute the mean ISE of both estimators which is given by $$\mathrm{ISE}\left(\widehat{F}\right):=\displaystyle\int_{-1}^1\left(\widehat{F}(x)-F(x)\right)^2d\,x,$$
					for each example of the distribution function and across various sample sizes $n$ and $N$. The obtained values are summarized in Table \ref{tab1}, where we also present the optimal $m$ values corresponding to each simulated case.
					\begin{table}[H]\label{ise}
						\centering
						\caption{Mean ISE of $\widehat{F}_{m,n,N}$ and $\widehat{F}_{h}$. The bold values indicate the smallest values of mean ISE.}
						\begin{tabular}{ccccc}
							\toprule
							\textbf{Distribution} & \textbf{ $(N, n)$} &\textbf{mean $m_{\text{opt}}$} &\textbf{Lagrange Estimator} & \textbf{Kernel Estimator} \\
							\midrule
							\multirow{3}{*}{$\mathcal{B}(2,2)$} & (100, 250)  & 7.2 & {\bf0.00271}  & 0.00772  \\
							& (250, 500) & 11.5 & {\bf0.00148}  & 0.00468  \\
							& (500, 1000) & 9.2 & {\bf0.00066}  & 0.00164  \\
							\midrule
							\multirow{3}{*}{$\mathcal{G}(2,1)$}  & (100, 250) & 9.5  & {\bf0.00941}  & 0.01601  \\
							& (250, 500)& 12.5  & {\bf 0.00509}  &  0.00697  \\
							& (500, 1000) & 14.5& {\bf0.00197}  & 0.00278  \\
							\midrule \multirow{3}{*}{$\mathcal{N}(0.5,0.25) $}
							& (100, 250) & 8.5 &  {\bf0.00393} & 0.01891  \\
							& (250, 500) & 10 & {\bf0.00111}  & 0.01036 \\
							& (500, 1000) & 11.5 & {\bf0.00095}  & 0.00957 \\
							\midrule \multirow{3}{*}{$0.5\mathcal{N}$(-2,1) + 0.5$\mathcal{N}(3,0.5)  $}
							& (100, 250) & 14 & {\bf0.00953}  & 0.05338 \\
							& (250, 500) & 14 & {\bf0.01000}  & 0.04976  \\
							& (500, 1000) & 14.5 & {\bf0.00101}  & 0.04592  \\
							\bottomrule
						\end{tabular}
						\label{tab1}
					\end{table}
					To illustrate the numerical results of Table~\ref{tab1}, a qualitative comparison of the mean ISE performance of both estimators is presented in the following figure.
					\begin{figure}[H]
						\centering
						\includegraphics[width=1\linewidth]{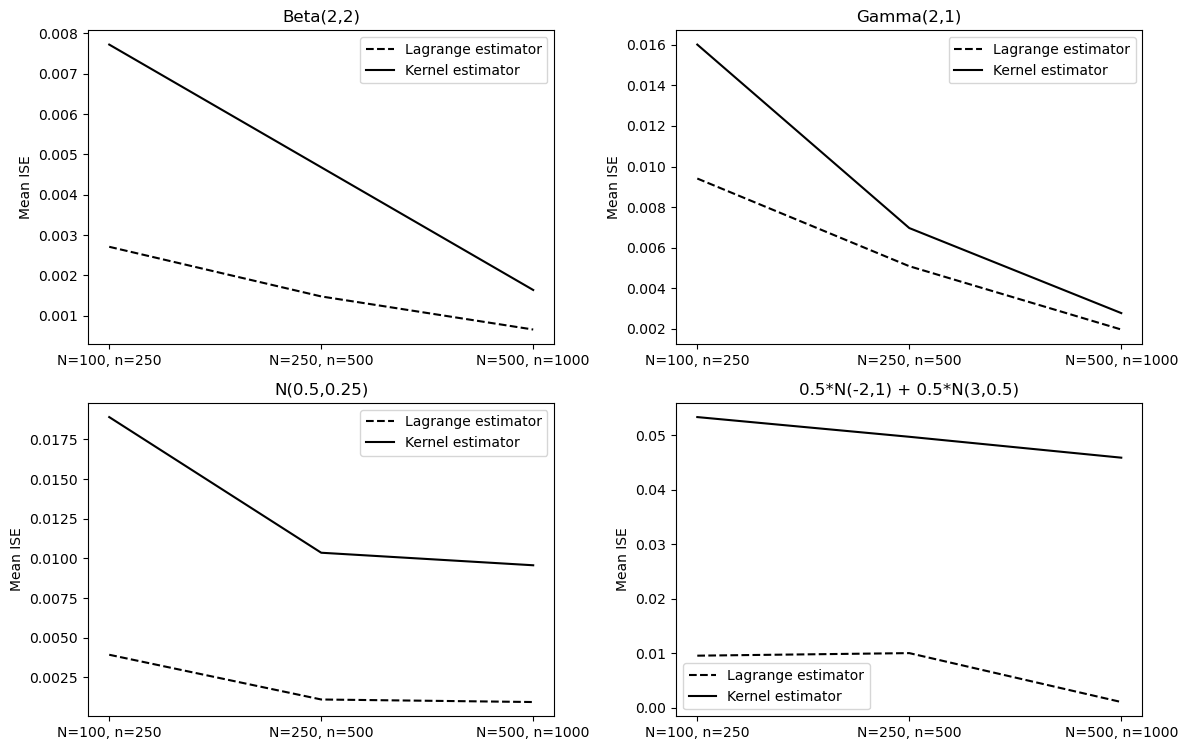}
						\caption{Qualitative comparison between the mean ISE of both estimators for each distribution function}
						\label{mean_ISE}
					\end{figure}

					\paragraph{Interpretations:}
					\begin{itemize}
						\item  From Figures \ref{fig1}-\ref{fig4}, we first observe that the accuracy of both estimators improves with a larger sample size. We also observe that the Lagrange estimator $\widehat{F}_{m,n,N}$ captures local curvature better and consistently follows the general shape of the true distribution function more smoothly than the kernel estimator $\widehat{F}_{h}$.
						\item The Lagrange estimator provides stable boundary behavior across all figures, while Kernel estimator shows more fluctuations and larger bias near the boundaries, particularly for the Beta and Gamma distributions in Figures \ref{fig1}-\ref{fig2}. For these examples, the Kernel estimator overestimates the true distribution function at the lower tail and underestimates the upper tail especially for small sample sizes.
						\item For the normal distribution plotted in Figure \ref{fig3}, it's clearly observable that the kernel estimator overestimates the true distribution near the bound $x=1$.
						\item As illustrated in Table \ref{tab1}, the mean optimal order $m$ increases as $n$ and $N$ increase.
						\item Across numerical results presented in Table \ref{tab1} and qualitative comparison illustrated by Figure \ref{mean_ISE}, the mean ISE decreases as the number of subjects $N$ and the observations per subject $n$ increase, indicating that both estimators improve with larger sample sizes. It's also obvious that the Lagrange estimator exhibits lower mean ISE.
						\item Figure \ref{mean_ISE} illustrates that, for all distributions examples and all sample sizes, the minimal mean ISE-value of the Lagrange estimator is always lower than that of the Kernel estimator.
						\item Across all results, the Lagrange estimator outperforms the kernel estimator.
						
						\item From Table \ref{tab2}, it is obvious that all estimators improve both precision and accuracy with a larger sample size ($N,n$). Specifically, 
						we can observe that for all parameters, the mean estimates are consistently close to their true values.
						\item By increasing $(N,n)$ from $(100,250)$ to $(500, 1000)$ the standard deviations of $\widehat{H}$, $\widehat{\gamma}^2$ and $\widehat{\sigma}^2$ decrease. This reflects higher precision and illustrates the theoretical asymptotic consistency of estimators.
						\item Across all distributions, the standard deviation of $\widehat{\phi}_i$ decreases as both the number of subjects $N$ and the number of observations per subject $n$ increase. This implies that the variance of $\widehat{\phi}_i$ tends to zero as the sample size increases.
						\item All the numerical results reported in Table \ref{tab2} confirm the theoretical asymptotic properties we investigated earlier, which validates the effectiveness of the proposed parametric estimation procedure.
					\end{itemize}
					\section{Empirical Application to Cryptocurrency Markets}\label{empiri}
					This section illustrates the empirical performance of our estimation methodology using cryptocurrency data. These assets present a perfect testing ground due to their several unique statistical characteristics. Indeed, they show strong cross-sectional heterogeneity resulting from fundamentally different blockchain technologies, have clear long-memory characteristics (as we will demonstrate), and have sharply defined volatility regimes, providing a rigorous environment to validate our mixed fractional Brownian motion model with random effects.
					We use daily historical prices from Yahoo Finance for five large cryptocurrencies (AVAX, ETH, LINK, LTC, and MATIC) from 2018-01-01 to 2025-01-01 and form log-returns \(r_t \;=\; \log\!\left(\frac{P_t}{P_{t-1}}\right)\), where $P_t$ is the daily adjusted closing price at time $t$. To estimate global parameters on stable blocks while preserving cross-sectional heterogeneity, we construct a pseudo-panel of $N\in\{100, 400, 900\}$ trajectories, each of length $n\in\{250, 500, 1000\}$ daily returns. Concretely, we set
					$r^{i}=(r^{i}_1,\ldots,r^{i}_n), k=1,\ldots,n,$ where each
					$r^{i}_k, \; i=1,\ldots,N$ corresponds to the increment $\Delta Y_{k}^i, i=1,\ldots,N $ in our theoretical model, 
					so each asset yields many short trajectories drawn from the same time series. The sampling step is fixed at \(h=\frac{1}{252}\). The global parameters \((H,\gamma^2,\sigma^2)\) are estimated using moments computed from \(r^{i}_k\) and applying our generalized method of moments developed earlier. Then, within the \(i\)-th trajectory, we estimate the random effects via
					\[
					\widehat{\theta}_i=\frac{1}{nh}\sum_{k=1}^n r^{i}_k, 
					\qquad 
					\widehat{\phi}_i=\widehat{\theta}_i+\tfrac{1}{2}\,\widehat{\sigma}^2.
					\] 
					Finally, the distribution of random effects is recovered, on the basis of $\widehat{\phi}_i, i=1,\ldots, N$, by the Lagrange estimator \eqref{F_phihat}. 
					
					
					The results of the estimation of the global parameters for each cryptocurrency are presented in the following table.

\begin{figure}[H]
						\centering
						\includegraphics[width=1\linewidth]{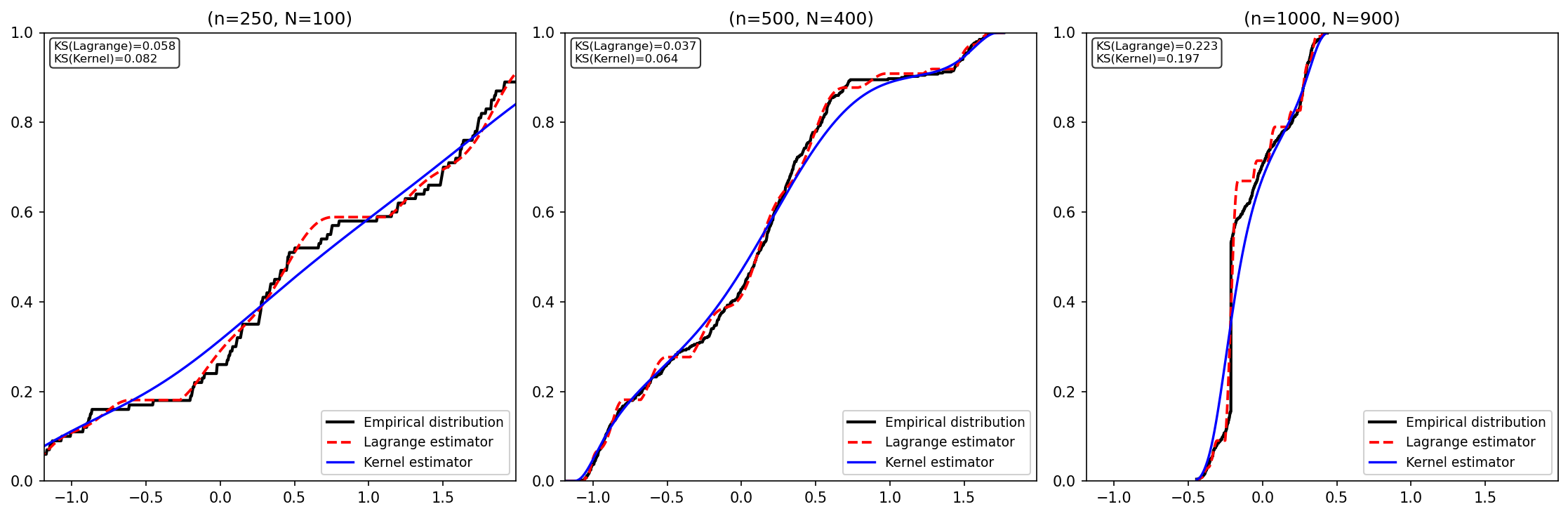}
						\caption{Qualitative comparison between random effects Distribution and its Lagrange and Kernel estimators across different sample sizes for the AVAX-USD.}
						\label{avax}
					\end{figure}

\begin{figure}[H]
						\centering
						\includegraphics[width=1\linewidth]{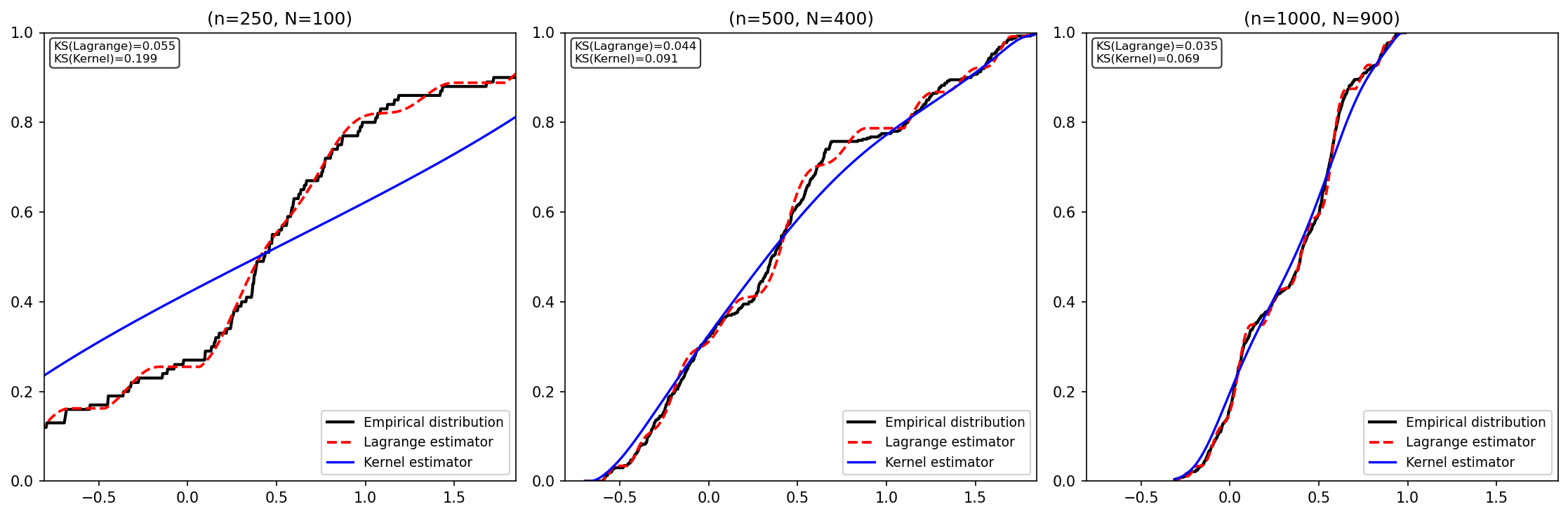}
						\caption{Qualitative comparison between random effects Distribution and its Lagrange and Kernel estimators across different sample sizes for the ETH-USD.}
						\label{eth}
					\end{figure}

\begin{figure}[H]
						\centering
						\includegraphics[width=1\linewidth]{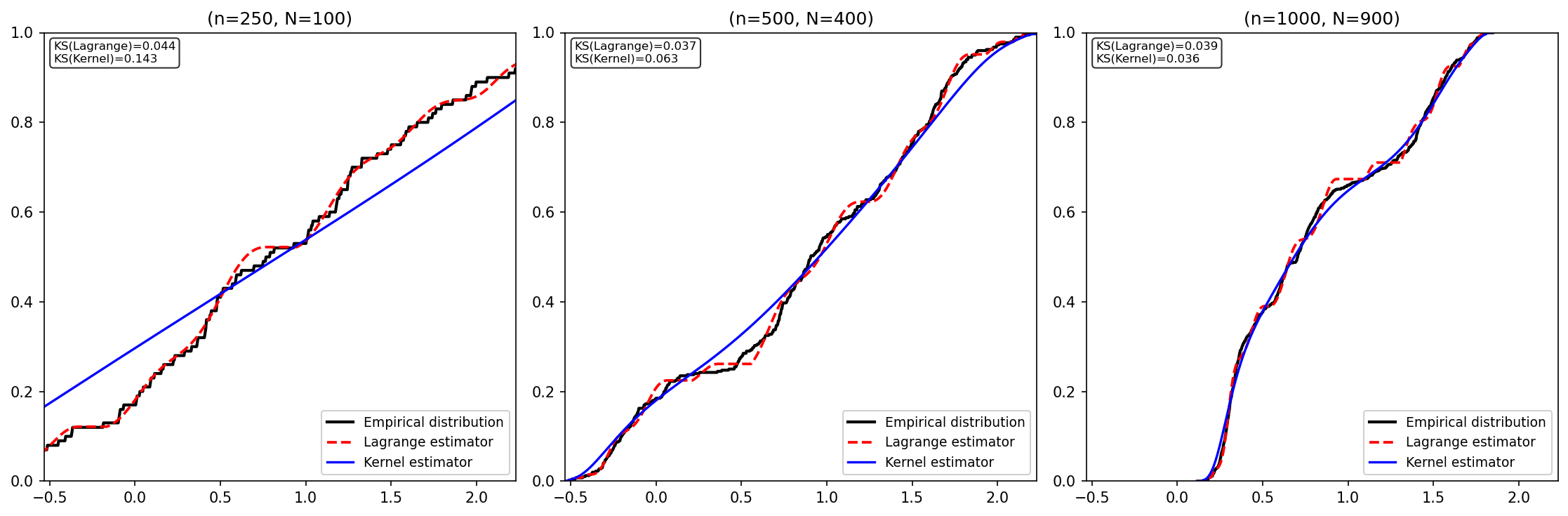}
						\caption{Qualitative comparison between random effects Distribution and its Lagrange and Kernel estimators across different sample sizes for the LINK-USD.}
						\label{link}
					\end{figure}

\begin{figure}[H]
						\centering
						\includegraphics[width=1\linewidth]{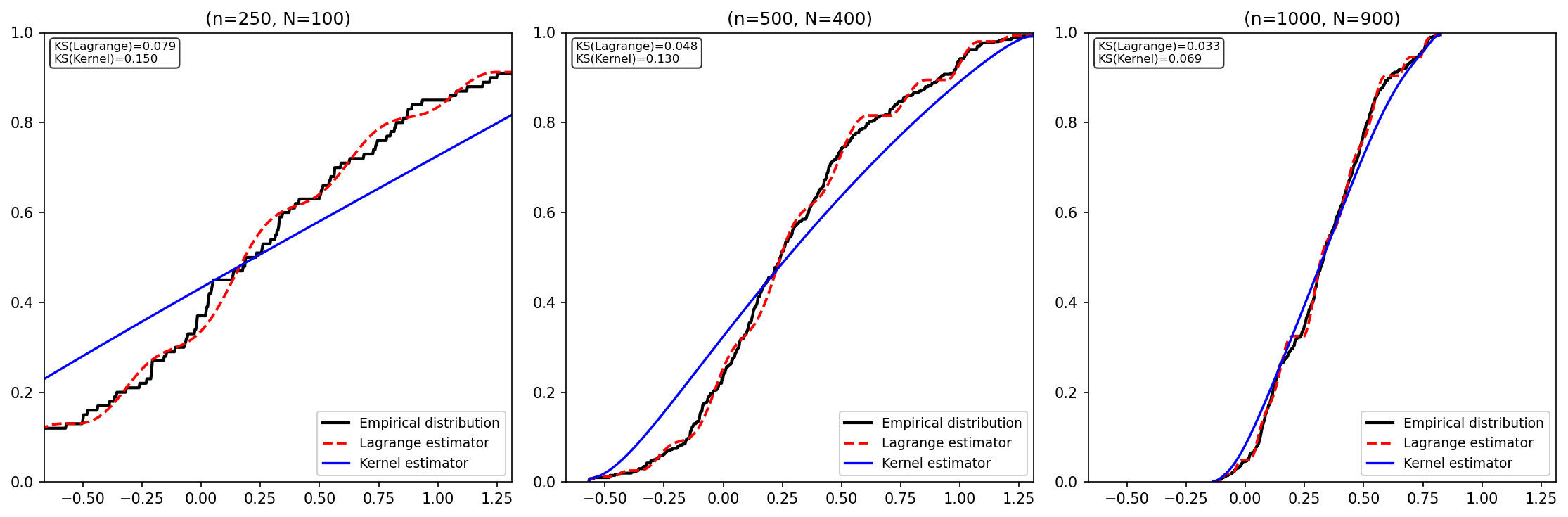}
						\caption{Qualitative comparison between random effects Distribution and its Lagrange and Kernel estimators across different sample sizes for the LTC-USD.}
						\label{ltc}
					\end{figure}

\begin{figure}[H]
						\centering
						\includegraphics[width=1\linewidth]{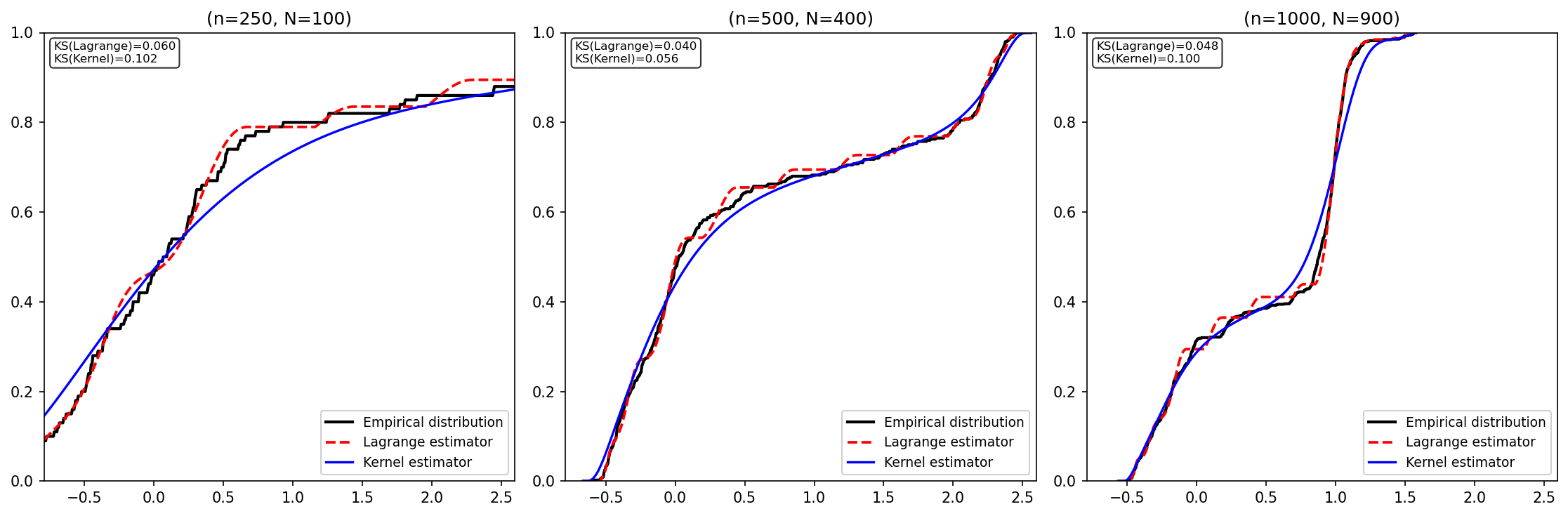}
						\caption{Qualitative comparison between random effects Distribution and its Lagrange and Kernel estimators across different sample sizes for the MATIC-USD.}
						\label{matic}
					\end{figure}

					\begin{table}[H]
						\centering
						\caption{Parameter estimates across different sample sizes}
						\label{tickers_param}
						\begin{tabular}{lrrrrrr}
							\toprule
							Ticker & $(n,N)$ & $\widehat{H}$ & $\widehat{\gamma}^2$ & $\widehat{\sigma}^2$ & $\widehat{\phi}_{\text{mean}}$ & $\widehat{\phi}_{\text{std}}$  \\
							\midrule
							& (250,100)  & 0.638877 & 0.000723 & 0.513261 & 0.442342 & 0.912931  \\
							ETH-USD  & (500,400)  & 0.639164 & 0.000649 & 0.519850 & 0.527801 & 0.633322  \\
							& (1000,900) & 0.638413 & 0.000667 & 0.554567 & 0.618960 & 0.310896  \\
							\midrule
							& (250,100)  & 0.604507 & 0.668701 & 0.370877 & 0.137566 & 0.759827  \\
							LTC-USD  & (500,400)  & 0.634520 & 0.610994 & 0.459273 & 0.226172 & 0.397326  \\
							& (1000,900) & 0.627225 & 0.680661 & 0.467827 & 0.249960 & 0.214846 \\
							\midrule
							& (250,100)  & 0.649911 & 0.000673 & 0.983277 & 0.698825 & 1.149438  \\
							AVAX-USD & (500,400)  & 0.610892 & 0.000915 & 0.895812 & 0.504652 & 0.701866  \\
							& (1000,900) & 0.609321 & 0.000769 & 0.791624 & 0.315646 & 0.218235  \\
							\midrule
							& (250,100)  & 0.633006 & 0.002370 & 1.222424 & 1.109158 & 1.462408  \\
							MATIC-USD & (500,400)  & 0.633170 & 0.002332 & 1.209355 & 1.183159 & 1.044514  \\
							& (1000,900) & 0.635014 & 0.002165 & 1.239305 & 1.164758 & 0.581617  \\
							\midrule
							& (250,100)  & 0.645225 & 0.000686 & 0.917443 & 0.843615 & 0.924283  \\
							LINK-USD & (500,400)  & 0.649799 & 0.000576 & 0.927584 & 0.894226 & 0.710103  \\
							& (1000,900) & 0.658765 & 0.000451 & 0.966385 & 0.832600 & 0.503922 \\
							\bottomrule
						\end{tabular}
					\end{table}
					As both the number of trajectories $N$ and the number of discrete $n$ increase, global parameters $(\widehat{H}, \widehat{\gamma}^2, \widehat{\sigma}^2)$ tend to stabilize their values across all cryptocurrencies, thus, confirming the theoretical consistency of the two-step asymptotic framework we introduced in section \ref{sec3_chapi3}.
					The estimated Hurst exponents cluster between 0.61 and 0.66 for most of assets, which indicates that cryptocurrency returns exhibit long-range dependence and match our theoretical assumption $H\in(\frac12,\frac34)$.  
					Hence, we can deduce that the shocks in the returns are carried forward over medium horizons rather than disappearing instantly, which is exactly a feature that the fractional component of the model captures.
					The fractional variance parameter $\widehat{\gamma}^2$ demonstrates a wide range of variability among the assets and is quite low for ETH, LINK, and AVAX and relatively high for MATIC and especially LTC, thus implying that some assets have strong memory. On the other hand, the short-term volatility parameter $\widehat{\sigma}^2$ is higher for MATIC and LINK and lower for both LTC and ETH.
					Moreover, the standard deviation of the random effects, $\widehat{\phi}_{std}$, gets considerably smaller as we use more data. This very important decrease shows that with a larger sample size, our estimation method distinguishes better real, lasting differences in the trend of an asset from random noise. In small samples, random effects reflect not only actual differences across trajectories but also estimation noise. In contrast, in larger samples, the noise component gets smaller and the remaining dispersion corresponds to genuine, persistent asset-specific drifts. \\
					Now, we present in Figures \ref{avax}–\ref{matic} a qualitative comparison between the empirical distribution functions of the random effects for each asset, along with their Lagrange estimator \eqref{F_phihat} and kernel estimator \eqref{kern_phi}, across three sample sizes $(n,N)=(250,100),(500,400),(1000,900)$.


					Across all assets, the Lagrange estimator (red dashed) is very close to the empirical distribution function (black line) even for relatively small samples, and its fit gets better in a very orderly manner as both $n$ and $N$ increase. This serves as a confirmation of the consistency and smoothness properties that were developed in our asymptotic theory: our proposed estimator is closer to the true distribution of the random effects when there is more data both within and across trajectories. On the other hand, the kernel estimator (blue curve) is generally more biased in small samples, especially at the tails of the distribution, where boundary effects cause the curvature and flattening to be quite pronounced. 
					The Kolmogorov–Smirnov (KS) distances presented in each subplot of Figures \ref{avax}–\ref{matic} decrease for both estimators with larger samples, which indicates higher accuracy and less estimation errors. Particularly, our proposed estimator exhibits smaller KS distances than the kernel method, across almost all cryptocurrencies, which ensures better performance and greater stability close to the support boundaries.\\
					In summary, all these empirical results demonstrate the robustness performance of our proposed estimation methods. The stabilization of global parameters with increasing sample size confirms theoretical consistency, the significant reduction in $\phi_{\text{std}}$ validates the precision of our random effects estimator and confirms that the mixed fractional Brownian motion with random effects successfully captures the persistent and heterogeneous nature of cryptocurrency returns.

					\section{conclusion}\label{conclusion_chapi3}
					In this work, we investigated a hybrid estimation framework for a mixed fractional Black-Scholes model incorporating both fixed global parameters and subject-specific random effects. 
					By using moment-based statistics and discrete observations, we developed consistent and asymptotically normal estimators for model fixed parameters and random effects.
					We then addressed the estimation of the random effects distribution, particularly in cases where its support is compact, and we introduced a nonparametric approach based on Lagrange interpolation at Chebyshev-Gauss nodes. The estimator was shown to perform reliably and consistently as both the number of trajectories $N$ and the number of discrete observations per-subject increase. The robustness of our procedure was rigorously validated through a simulation study and an empirical application to cryptocurrency markets, demonstrating its practical utility in decomposing complex financial time series into persistent noise, transient volatility, and asset-specific fundamental trends. Extensions of the present work could address the more challenging case $H\geq \frac 34$, where the asymptotic limits become non‑Gaussian, requiring different normalizations and inference techniques.
					
					\section{Proofs of main results}
					\label{sec6_chapi3}
					\subsection{Proof of Lemma \ref{lem_ergod}}
					For each fixed $i = 1, \ldots, N$, $\theta_i$ is constant within each subject. Furthermore, the sequence $\{\Delta Y_k^i - \theta_i h, 0 \leq k \leq n-1\}$ can be written as $\{\Delta M_k^{H,i}, 0 \leq k \leq n-1\}$, where $\Delta M_k^{H,i} = \sigma \Delta B_k^i + \gamma \Delta B^{H,i}_k$. 
					As we have $\{\Delta M_k^{H,i}, 0 \leq k \leq n-1\}$ is a stationary Gaussian sequence, it remains only to verify the ergodicity condition. The autocovariance function is given by
					\begin{eqnarray}\label{ro}
						\rho(k) = \mathrm{Cov}(\Delta M^{H,i}_k, \Delta M^{H,i}_0) = 
						\begin{cases}
							\sigma^2 h + \gamma^2 h^{2H} & \text{si } k = 0, \\
							\gamma^2 \rho_H(k) & \text{si } k \neq 0,
						\end{cases}
					\end{eqnarray}
					where $\rho_H(k) = \frac{1}{2}\left(|k+1|^{2H} + |k-1|^{2H} - 2|k|^{2H}\right)$ for $k \geq 1$.
				In the case $H \in \left(\frac{1}{2}, 1\right)$, the asymptotic behavior of $\rho_H(k)$ is as follows
				\begin{eqnarray}\rho_H(k) \sim h^{2H} H(2H-1) k^{2H-2} \; \mbox{ as } \; k \to \infty.\label{simrho}\end{eqnarray}
				Thus
				$$\sum_{k=-\infty}^{\infty} |\rho(k)| = \sigma^2 h + 2\gamma^2 \sum_{k=1}^{\infty} |\rho_H(k)| < \infty,$$
				since $\sum\limits_{k=1}^{\infty} k^{2H-2} < \infty$ for $H \in \left(\frac{1}{2}, 1\right)$ because $2H-2 > -1$.
				Therefore $\rho(k) \to 0$ as $k \to \infty$ and $\sum|\rho(k)| < \infty$. By Corollary \ref{corA1} in the appendix,  the sequence $\{\Delta M_k^{H,i}, 0 \leq k \leq n-1\}$ is ergodic.

				\subsection{Proof of Proposition \ref{prop1_chap3}}
				For each $i=1,\ldots,N$, $\widehat{\theta}_i= \dfrac{1}{nh}\sum\limits_{k=0}^{n-1}\Delta Y^i_k=\theta_i + \dfrac{1}{nh}\sum\limits_{k=0}^{n-1}\Delta M^{H,i}_k$. 
				By Lemma~\ref{lem_ergod}, the Gaussian sequence $\{\Delta M^{H,i}_k,\ 0\le k\le n-1\}$ is stationary and ergodic. Hence, by Theorem \ref{erg_th}, 
				\[
				\frac{1}{n}\sum_{k=0}^{n-1}\Delta Y^i_k
				\xrightarrow{a.s.} \theta_i\, h, \quad \mathrm{as}\quad n\to\infty.
				\]
				which proves that
				$
				\widehat{\theta}_i
				\;\xrightarrow{a.s.}\; \theta_i,$ as $n\to\infty$.\\
				We now study the asymptotic normality of $\widehat{\theta}_i$. 
				For each $i=1,\ldots, N$, we have
				\[
				\widehat{\theta}_i-\theta_i \;=\; \frac{\sigma}{nh}\sum_{k=0}^{n-1}\Delta B_k^i \;+\; \frac{\gamma}{nh}\sum_{k=0}^{n-1}\Delta B^{H,i}_k.
				\]
				For the Brownian part, since the increments of $B^i$ are independent and stationary, we obtain
				\[
				\Var\!\left(\frac{\sigma}{nh}\sum_{k=0}^{n-1}\Delta B^i_k\right)
				= \frac{\sigma^2}{(nh)^2}\Var\!\left(B^i_{nh}-B^i_0\right)
				= \frac{\sigma^2}{(nh)^2}nh
				= \sigma^2 (nh)^{-1},
				\]
				which implies that the Brownian motion term is of order $n^{-\frac12}$.\\
				For the fBm term, we have
				\[
				\Var\!\left(\frac{\gamma}{nh}\sum_{k=0}^{n-1}\Delta B^{H,i}_k\right)
				= \frac{\gamma^2}{(nh)^2}\Var\,\big(B^{H,i}_{nh}-B^{H,i}_0\big)
				= \frac{\gamma^2}{(nh)^2}\,(nh)^{2H}
				= \gamma^2 h^{2H-2} n^{2H-2}.
				\]
				Thus, the fBm term behaves like $n^{H-1}$. Since $H>\tfrac12$, we obtain
				\[
				\frac{n^{H-1}}{n^{-\frac12}} = n^{H-\tfrac12}\;\longrightarrow\;\infty,
				\]
				which shows that the fBm term asymptotically dominates the Brownian one. Therefore, the effective normalization for $\widehat{\theta}_i-\theta_i$ is $n^{1-H}$.\\
				Now, multiplying the error term by $n^{1-H}$, we obtain
				\begin{equation}
					n^{1-H}\big(\widehat{\theta}_i-\theta_i\big)
					= \frac{\sigma}{h}\,n^{-H}\big(B^i_{nh}-B^i_0\big)
					+ \frac{\gamma}{h}\,n^{-H}\big(B^{H,i}_{nh}-B^{H,i}_0\big).\label{theta_err}
				\end{equation}
				For the first term of \eqref{theta_err}, \(
				\Var\!\left(\frac{\sigma}{h}n^{-H}\big(B^i_{nh}-B^i_0\big)\right)
				= \frac{\sigma^2}{h^2}\,n^{-2H} nh
				= \frac{\sigma^2}{h}\,n^{1-2H},
				\)
				and since $H>\frac12$, we obtain \begin{eqnarray}\label{p_bm}\frac{\sigma}{h}n^{-H}(B^i_{nh}-B^i_0)\;\xrightarrow{\mathbb{P}}\;0\quad \mathrm{as}\quad n\to\infty.\end{eqnarray}
				For the second term of \eqref{theta_err},
				since $B^H_{nh}-B^H_0\sim \mathcal{N}(0,(nh)^{2H})$, we obtain
				\begin{eqnarray}
					\frac{\gamma}{h}n^{-H}\big(B^{H,i}_{nh}-B^{H,i}_0)\;\xrightarrow{d}\;\mathcal{N}\big(0, \gamma^2 h^{2H-2}\big).\label{d_fbm}\end{eqnarray}
				Replacing \eqref{p_bm} and \eqref{d_fbm} into \eqref{theta_err} and using Slutsky’s theorem, we obtain
				\[
				n^{1-H}\big(\widehat{\theta}_i-\theta_i\big)\;\xrightarrow{d}\;\mathcal{N}\!\big(0,\gamma^2 h^{2H-2}\big),
				\]
				which completes the proof.

				\subsection{Proof of Proposition \ref{cvg.as}}
				The idea of this proof is first to establish the almost sure convergence of the per-subject statistics $\xi_{n}^i$, $\eta_{n}^i$ and $\zeta_{n}^i$ when the number of observations per subject increases to infinity by using the ergodic theorem \ref{erg_th} and then deduce the result of Proposition \ref{cvg.as} by using the SLLN applied to the empirical statistics $\bar{V}_{n,N}$, $\bar{\xi}_{n,N}$, $\bar{\eta}_{n,N}$ and $\bar{\zeta}_{n,N}$.\\ 
				Since each per-subject statistic $\xi^i_n$, $\eta^i_n$, $\zeta^i_n$ is constructed from the trajectory of subject $i$, and all subjects $\{Y^i, \ldots, Y^N\}$ are assumed i.i.d, it follows that, for fixed $n$, the sequences $\{\widehat\theta_i,\ i=1,\ldots, N\}$, $\{\xi^i_n,\ i=1,\ldots, N\}$, $\{\eta^i_n,\ i=1,\ldots, N\}$ and $\{\zeta^i_n,\ i=1,\ldots, N\}$ are i.i.d. In particular, the sequence $\{\widehat{\theta}^2_i,\ i=1,\ldots, N\}$ is also i.i.d., since it is obtained as a measurable transformation of the i.i.d. sequence $\{\widehat\theta_i,\ i=1,\ldots, N\}$. Now, in order to establish the almost sure convergence of these i.i.d. sequences, we apply the ergodic theorem to each per-subject statistic. 
				For each fixed $i=1,\ldots,N$, we rewrite the per-subject statistics as follows 
				\begin{eqnarray*}\xi_{n}^i&=&\dfrac{1}{n}\sum\limits_{k=0}^{n-1} \left[ \left(\Delta Y^i_k-\theta_i h\right)+\theta_i h\right]^2,\quad
					\eta_{n}^i=\dfrac{1}{n}\sum\limits_{k=0}^{n-1} \left[ \left(\Delta Y^i_k-\theta_i h\right)+\theta_i h\right]\left[ \left(\Delta Y^i_{k+1}-\theta_i h\right)+\theta_i h\right],\\
					\zeta_{n}^i&=&\frac1n\sum\limits_{k=0}^{n-1} \left[ \left(\Delta Y^i_k-\theta_i h\right)+ \left(\Delta Y^i_{k+1}-\theta_i h\right)+2\theta_i h\right]\left[ \left(\Delta Y^i_{k+2}-\theta_i h\right)+ \left(\Delta Y^i_{k+3}-\theta_i h\right)+2\theta_i h\right].
				\end{eqnarray*}
				By Lemma \ref{lem_ergod}, the centered sequence $\{\Delta Y^i_k-\theta_i h,\ 0\le k\leq n-1\}$ is stationary and ergodic. Then, by applying Theorem \ref{erg_th} conditionally on $\theta_i$, we obtain for each fixed $i=1,\ldots,N$ and as $n$ tends to infinity, 
				\begin{eqnarray*}
					\xi_{n}^i& \xrightarrow{\text{a.s.}}& \mathbb{E}\left[(\Delta Y_0^i)^2|\theta_i\right] =\theta_i^2 h^2+\sigma^2 h + \gamma^2 h^{2H}:=\xi_\infty^i,\\
					\eta_{n}^i&\xrightarrow{\text{a.s.}}& \mathbb{E}\left[\Delta Y_0^i \Delta Y_1^i|\theta_i\right]=\theta_i^2 h^2 +  \gamma^2 h^{2H} (2^{2H - 1} - 1):=\eta_\infty^i,\\
					\zeta_{n}^i &\xrightarrow{\text{a.s.}}& \mathbb{E}\left[(Y_{2h}^i - Y_0^i)(Y_{4h}^i - Y_{2h}^i)|\theta_i\right]=4 \theta_i^2 h^2 + \gamma^2 h^{2H}2^{2H} \left(2^{2H -1}-1 \right):=\zeta_\infty^i.
				\end{eqnarray*}
				Furthermore, since for each fixed \(i=1,\ldots, N\), $\widehat\theta_i \xrightarrow{a.s.} \theta_i$ as $n\to\infty$, by the continuous mapping theorem, we obtain $\widehat\theta_i^2 \xrightarrow{a.s.} \theta_i^2$ as $n\to\infty$.
				Denote by $U_\infty^i$ the almost sure limit vector of per-subject statistics for subject  $i$, i.e,
				\begin{eqnarray}
					U_\infty^i :=\Big(V_\infty^i, \xi_\infty^i, \eta_\infty^i, \zeta_\infty^i \Big)^\top.\label{u_i_inf}
					\end{eqnarray}
					
					We now turn to the almost sure convergence of the averaged statistics. The vectors \(U_\infty^i\) are i.i.d.\ across $N$ subjects and under assumption $(A1)$, they have finite means. Hence, the SLLN yields
					\[
					\bar U_{n,N} = \frac{1}{N}\sum_{i=1}^N U_\infty^i \xrightarrow{\text{a.s.}} \E[U_\infty^1] =: U_\infty,
					\]
					where \(U_\infty:=\Big(V_\infty, \xi_\infty, \eta_\infty, \zeta_\infty \Big)\), with 
					$V_\infty = \E[\theta_1^2],$ \quad
					$\xi_\infty = \E[\theta_1^2] h^2 + \sigma^2 h + \gamma^2 h^{2H},$\\
					$\eta_\infty = \E[\theta_1^2] h^2 + \gamma^2 h^{2H}(2^{2H-1}-1),$ \quad
					$\zeta_\infty = 4\E[\theta_1^2] h^2 + \gamma^2 h^{2H} 2^{2H}(2^{2H-1}-1).$ 
					This completes the proof.


					\subsection{Proof of Proposition \ref{prop3_chap3}}
					The proof of Proposition \ref{prop3_chap3} is structured in two steps. First, by applying the ergodic theorem \ref{erg_th} conditionally on $\theta_i$, we establish the almost sure convergence of the per-subject statistics as $n\to\infty$. Second, we derive the asymptotic normality of the averaged statistics by invoking the classical multivariate Central Limit Theorem (CLT) for i.i.d. vectors $U_\infty^i$, $i=1,\ldots,N$ as $N\to\infty$, and combine the two steps using Slutsky’s theorem.
					\paragraph{Asymptotic normality of $U^i_n:=\left( \widehat{\theta}_i^2,\, \xi_{n}^i, \eta_{n}^i, \zeta^i_n \right)^\intercal$}\ \\
					Let us define the following conditionally centered statistics
					\[
					\widetilde{\xi}_n^i := \xi_{n}^i - \mathbb{E}[\xi_{n}^i \mid \theta_i], \quad 
					\widetilde{\eta}_n^i := \eta_{n}^i - \mathbb{E}[\eta_{n}^i \mid \theta_i] \quad \mathrm{and}\quad
					\widetilde{\zeta}_n^i := \zeta_{n}^i - \mathbb{E}[\zeta_{n}^i \mid \theta_i].
					\]
					For each fixed $i$, these statistics can be decomposed into a linear term and a centered quadratic term as follows:
					\begin{eqnarray*}
						\widehat{\theta}_i^2-\theta_i^2 &= L_{n,\theta}+ Q_{n,\theta}&:=\dfrac{2\theta h}{n}\sum\limits_{k=0}^{n-1}\Delta M_k^{H,i}+ \Big(\dfrac{1}{nh}\sum_{k=0}^{n-1}\Delta M_k^{H,i}\Big)^2,\\
						\widetilde{\xi}_n^i&= L_{n,\xi}+ Q_{n,\xi}&:=\dfrac{2\theta h}{n}\sum\limits_{k=0}^{n-1}\Delta M_k^{H,i}+ \dfrac{1}{n}\sum_{k=0}^{n-1}\Big((\Delta M_k^{H,i})^2-\mathbb{E}[(\Delta M_k^{H,i})^2]\Big),\\
						\widetilde{\eta}_n^i&= L_{n,\eta}+ Q_{n,\eta}&:= \frac{\theta h}{n}\sum_{k=0}^{n-1}\big(\Delta M_k^{H,i}+\Delta M_{k+1}^{H,i}\big)
						+ Q_{n,\eta},\\
						\widetilde{\zeta}_n^i
						&= L_{n,\zeta}+ Q_{n,\zeta}&:= \frac{2\theta h}{n}\sum_{k=0}^{n-1}\big(\Delta M_k^{H,i}+\Delta M_{k+1}^{H,i}+\Delta M_{k+2}^{H,i}+\Delta M_{k+3}^{H,i}\big)
						+ Q_{n,\zeta}.
					\end{eqnarray*} 
					For each statistic, the linear term is of Hermite rank 1, while the centered quadratic term is of Hermite rank 2. So, to study the asymptotic distribution of the linear term, we will use the Gaussian properties of the sequence $\{\Delta\, M_k^{H, i},\ 0\le k\le n-1\}$. Then, we will apply the Breuer-Major theorem to establish the asymptotic normality of the centered quadratic term.
					
					As proved previously, the appropriate normalization for the linear term involving $\Delta M^{H,i}_k$ is $n^{\,1-H}$, since the fractional Brownian component dominates the asymptotic behavior. Then, by scaling the linear terms and using the fact that $\dfrac{1}{n} \sum\limits_{k=0}^{n-1}\Delta M_k^{H, i}\ \xrightarrow{d}\ \mathcal{N}\Big(0, \gamma^2\,h^{2H}\,n^{2H-2}\Big)$, we obtain 
					\begin{align*}
						n^{\,1-H} L_{n,\theta}&\ \xrightarrow{d}\ \mathcal{N}\Big(0, 4\theta_i^2\gamma^2\,h^{2H-2}\Big),\\
						n^{1-H} L_{n,\xi}&\ \xrightarrow{d}\ \mathcal{N}\Big(0, 4\theta_i^2\gamma^2\,h^{2H+2}\Big),\\
						n^{1-H} L_{n,\eta}&\ \xrightarrow{d}\ \mathcal N\big(0,\ 4\theta_i^2\gamma^2 h^{2H+2}\big),\\
						n^{1-H} L_{n,\zeta} &\xrightarrow{d}\mathcal N\big(0,\ 64\theta_i^2\gamma^2 h^{2H+2}\big).
					\end{align*}
					
					We now proceed to study the quadratic terms. For each per-subject statistic, the quadratic term \(Q_{n,}\) is a centered quadratic functional of the stationary Gaussian sequence \(\{\Delta M_k^{H,i}, \ 0\le k\le n-1\}\), which implies that each quadratic term is of Hermite rank 2. Then, to apply the Breuer-Major theorem, we need to verify that $\sum\limits_{k \in \mathbb{Z}} |\rho(k)|^2 < \infty$. 
					Using the expression \eqref{ro}, we obtain
					\begin{align*}
						\sum_{k \in \mathbb{Z}} |\rho(k)|^2 &= |\sigma^2 h|^2 + \sum_{k \neq 0} |\gamma^2 \rho_H(k)|^2\\
						&= \sigma^4 h^2 + \gamma^4 \sum_{k \neq 0} |\rho_H(k)|^2.
					\end{align*}
					From the asymptotic behavior of $\rho_H(k)$ given by \eqref{simrho}, we deduce
					$$\sum_{k \neq 0} |\rho_H(k)|^2 \sim C h^{4H} \sum_{k \neq 0} |k|^{4H-4},$$
					for some constant $C>0$. This series converges if and only if $4H - 4 < -1$, which gives us $H < \frac{3}{4}$. Hence, for $H\in\big(\frac12, \frac34\big)$, the Breuer-Major theorem applies for each quadratic term and yields
					$$\sqrt{n}\,Q_{n,\cdot} \xrightarrow{d} \mathcal N(0,v_Q).$$
					Scaling by $n^{1-H}$ and since $n^{\frac12-H}\to0$ for $H>\frac12$, we obtain 
					\[
					n^{1-H}\,Q_{n,\cdot} = \frac{n^{1-H}}{\sqrt{n}}\,\sqrt{n}\,Q_{n,\cdot} = \sqrt{n}\,Q_{n,\cdot} n^{\frac12-H} \overset{\mathbb{P}}{\longrightarrow} 0,\quad \mathrm{as}\quad n\to\infty.
					\]
					Therefore, the quadratic terms vanish under the $n^{1-H}$-normalization. 
					Now, to prove the joint asymptotic normality of all statistics,  we consider, for any \( \alpha = (\alpha_1, \alpha_2, \alpha_3, \alpha_4)^\intercal \in \mathbb{R}^4 \), the linear combination \[
					\mathcal{L}_n:=\alpha_1\left(\widehat{\theta}^2_i-\theta^2_i\right) + \alpha_2 \widetilde{\xi}_n^i + \alpha_3 \widetilde{\eta}_n^i + \alpha_4 \widetilde{\zeta}_n^i .
					\]
					After scaling by $n^{1-H}$ and using the previous results, we obtain
					\[
					n^{1-H}\ \alpha^\top \mathcal{L}_n\ \overset{d}{\longrightarrow}\ \mathcal{N}\left(0,\ \alpha^\top \Sigma_{\theta_i}\,\alpha \right),\quad \text{as}\quad n\to\infty,
					\]
					where $\Sigma_{\theta_i}:= v(\theta_i)v(\theta_i)^\top$ with $v(\theta_i):=\Big(4\theta_i^2\gamma^2\,h^{2H-2},\ 4\theta_i^2\gamma^2\,h^{2H+2},\ 4\theta_i^2\gamma^2\,h^{2H+2},\ 64\theta_i^2\gamma^2\,h^{2H+2}\Big)^\top$ is the vector of variances of the linear terms.
					Therefore, by the Cramèr-Wold device, the joint asymptotic normality of the vector $U_{n}^i=\Big(\widehat{\theta}_i^2, \xi_{n}^i, \eta_{n}^i, \zeta_{n}^i\Big)^\top$ holds, that is,\begin{eqnarray}
						n^{1-H}
						\big(U_{n}^i-\E[U_{n}^i\mid\theta_i]\big)\xrightarrow{d}\mathcal N\big(0,\Sigma_{\theta_i}\big),\quad \mathrm{as}\quad n\to\infty\label{asymp_i}.
					\end{eqnarray}
					\begin{rem}
						The normalization factor $n^{1-H}$ in the equation \eqref{asymp_i} (rather than the standard $\sqrt{n}$) arises because when $H > \frac{1}{2}$, the fBm component, which contributes errors of order $n^{H-1}$, dominates the Brownian motion component, which contributes errors of order $n^{-1/2}$. This different scaling is characteristic of statistical inference for fractional processes and distinguishes it from classical Brownian-based inference.
					\end{rem}

				\paragraph{Asymptotic normality of $\bar{U}_{n,N}:=\left( \bar{V}_{n,N},\ \bar{\xi}_{n,N},\ \bar{\eta}_{n,N},\ \bar{\zeta}_{n,N}\right)^\intercal$}\ \\
				For each fixed $i$, we decompose the vector of per-subject statistics as follows $$U_{n}^i=\mathbb{E}\left(U_{n}^i|\theta_i\right)+ \Big(U_{n}^i-\mathbb{E}\left(U_{n}^i|\theta_i\right)\Big).$$
				Averaging over $N$, we obtain
				\begin{eqnarray}\sqrt{N}\left(\bar{U}_{n,N}-U_\infty\right)=\dfrac{1}{\sqrt{N}} \sum\limits_{i=1}^N\Big(\mathbb{E}\left(U_{n}^i|\theta_i\right)-U_\infty\Big)+ \dfrac{1}{\sqrt{N}} \sum\limits_{i=1}^N\Big(U_{n}^i-\mathbb{E}\left(U_{n}^i|\theta_i\right)\Big).\label{sbar}\end{eqnarray}
				From the joint asymptotic normality of the per-subject statistics established in \eqref{asymp_i}, we deduce that  
				\[
				U_{n}^i - \E[U_{n}^i\mid\theta_i] = O_{\mathbb{P}}\big(n^{H-1}\big),
				\]
				and since for $H\in\big(\frac12,\frac34\big)$, \(n^{H-1}\to 0\) as $n\to\infty$, it follows that \begin{eqnarray}
					U_{n}^i - \E[U_{n}^i\mid\theta_i] = o_{\mathbb{P}}\big(1\big)\quad \mathrm{as}\ n\to\infty.\label{asymp_N1}
				\end{eqnarray}
				Therefore, as $n\to\infty$ for each fixed \(N\), the second term on the right-hand side of \eqref{sbar} vanishes. 
				\\ 
				Now, we proceed to study the asymptotic behavior of the first term on the right-hand side of \eqref{sbar}. For each fixed $i$, by the conditional ergodic theorem, we obtain
				\[
				\mathbb{E}[U_{n}^i \mid \theta_i] \xrightarrow{a.s.} U_\infty^i, \quad \text{as } n \to \infty,
				\]
				where $U_\infty^i$ is the per-subject limit vector defined in \eqref{u_i_inf}. Averaging over $N$ subjects, we obtain
				$$\dfrac{1}{\sqrt{N}} \sum\limits_{i=1}^N \Big(\mathbb{E}[U_{n}^i \mid \theta_i]- U_\infty\Big)\ \xrightarrow{a.s}\ \dfrac{1}{\sqrt{N}} \sum\limits_{i=1}^N \Big(U_\infty^i-U_\infty\Big),\quad \mathrm{as}\ n\to\infty.$$
				Furthermore, the vectors $U_\infty^i$, $i=1,\ldots,N$ are i.i.d, and under assumption $(A1)$ they have finite mean $U_\infty$. Then, the multivariate CLT yields
				\begin{eqnarray}
					\dfrac{1}{\sqrt{N}} \sum\limits_{i=1}^N \Big(U_\infty^i-U_\infty\Big)\ \xrightarrow{d}\ \mathcal{N}\big(0, \Sigma\big),\quad \mathrm{as}\ N\to\infty,\label{asymp_N2}
				\end{eqnarray}
				where $\Sigma:=\Cov\big(U^i_\infty\big)= \Var\big(\theta_i^2\big)\, \mathcal{H}\, \mathcal{H}^\top$, with $\mathcal{H}:=\big(1, h^2, h^2, 4h^2\big)^\top$.
				Finally, by combining \eqref{asymp_N1} and \eqref{asymp_N2}, the Slutsky's theorem \ref{slut} yields
				\begin{eqnarray}
					\sqrt{N}\left(\bar{U}_{n,N}-U_\infty\right)\ \xrightarrow{d} \mathcal{N}(0, \Sigma),\quad \mathrm{as}\quad n\to \infty \ \ \mathrm{then}\ \ \ N\to\infty.\label{asym_3}
				\end{eqnarray}
				
				\subsection{Proof of Proposition \ref{prop4_chap3}}
				To prove the results of this proposition, we begin with justifying the construction of global parameters, then we establish their strong consistency. 
				The estimators vector $\widehat\Theta$ is constructed from the limits \eqref{eq1}-\eqref{eq4} established in Proposition \ref{cvg.as}. In particular, from equations \eqref{eq3} and \eqref{eq4}, we obtain
				\[
				\eta_\infty - h^2 V_\infty \;=\; \gamma^2 h^{2H}\big(2^{2H-1}-1\big), 
				\qquad 
				\zeta_\infty - 4h^2 V_\infty \;=\; \gamma^2 h^{2H}\,2^{2H}\big(2^{2H-1}-1\big).
				\]
				Taking the ratio, we get
				$
				\frac{\zeta_\infty - 4h^2 V_\infty}{\eta_\infty - h^2 V_\infty} \;=\; 2^{2H},
				$ which implies that
				\[
				H = \frac{1}{2}\,\log_2^+\!\left(\frac{\zeta_\infty - 4h^2 V_\infty}{\eta_\infty - h^2 V_\infty}\right).
				\] 
				Replacing this quantity in \eqref{eq3}, we can express $\gamma^2$ as follows
				\[
				\gamma^2 \;=\; \frac{\eta_\infty - h^2 V_\infty}{h^{2H}\left(2^{2H-1}-1\right)}.
				\]
				Finally, substituting $H$ and $\gamma^2$ in \eqref{eq2}, we obtain 
				\[
				\sigma^2 \;=\; \frac{\xi_\infty - h^2 V_\infty - \gamma^2 h^{2H}}{h}.
				\]
				Now, since the quantities $V_\infty$, $\xi_\infty$, $\eta_\infty$ and $\zeta_\infty$ depend on the unknown expectation $\mathbb{E}\big(\theta_i^2\big)$, we replace them with their averaged statistics $\bar V_{n,N}, \bar\xi_{n,N}, \bar\eta_{n,N}, \bar\zeta_{n,N}$, which are based on the observed subjects. Therefore, we obtain the estimators $\widehat{H}$, $\widehat{\gamma}^2$ and $\sigma^2$ given by \eqref{hhat}, \eqref{gam} and \eqref{sig}.\\
				Now, we proceed to study the consistency of $\widehat\Theta$ which is based on results of Proposition \ref{cvg.as} combined with the continuous mapping theorem. 
				We consider the transformation $g$ defined as follows
				\begin{eqnarray*}
					g : \mathbb{R}^4 &\to& \mathbb{R}^3, \\
					x = \begin{pmatrix} x_1,\ x_2,\ x_3,\ x_4\end{pmatrix}^\top &\mapsto& 
					\begin{pmatrix}
						g_1(x),\
						g_2(x),\
						g_3(x)
					\end{pmatrix}^\top,
			\end{eqnarray*}
			where
			\[
			\begin{aligned}
				g_1(x) &= \dfrac{1}{2} \log_2^+ \left( \dfrac{x_4 - 4h^2 x_1}{x_3 - h^2 x_1} \right), \\
				g_2(x) &= \dfrac{x_3 - h^2 x_1}{h^{2g_1(x)} \left(2^{2g_1(x) - 1} - 1\right)}, \\
				g_3(x) &= \dfrac{x_2 - h^2 x_1 - g_2(x)  h^{2g_1(x)}}{h}.
			\end{aligned}
			\]
			Then, the estimator vector $\widehat{\Theta}$ can be expressed as 
			\(
			\widehat{\Theta} = g\big(\bar V_{n,N}, \bar\xi_{n,N},\bar\eta_{n,N},\bar\zeta_{n,N}\big)=g\left(\bar{U}_{n,N}\right).\) 
			To apply the continuous mapping theorem, we should first verify that the function $g$ is continuous in a neighborhood of the limit vector $U_\infty=\big( V_\infty, \xi_\infty,\eta_\infty, \zeta_\infty\big)$. 
			Examining the expressions for $g_1$, $g_2$, and $g_3$, the possible points of singularity are the denominators $x_3 - h^2 x_1$ and $2^{2g_1(x)-1}-1$ or $h^{2g_1(x)}$. Evaluating these quantities at $U_\infty$, we obtain
			\[
			\eta_\infty - h^2 V_\infty = \gamma^2 h^{2H}\big(2^{2H-1}-1\big).
			\]
			Since \(H>\tfrac12\) and \(\gamma>0\), it follows that \(\gamma^2 h^{2H}\big(2^{2H-1}-1\big)>0\), which ensures that \(x_3-h^2 x_1\) is nonzero at \(U_\infty\). Furthermore \(h^{2g_1(\mu)}=h^{2H}>0\) and \(2^{2H-1}-1>0\). Hence, \(g\) is well-defined and continuous in a neighborhood of \(U_\infty\). Therefore \(g\) is continuous at \(U_\infty\).\\  
			From Proposition \ref{cvg.as}, we have 
			\(
			\bar{U}_{n,N}\xrightarrow{a.s} U_\infty,
			\) as $n\to\infty$ first and then $N\to\infty$. 
			Then, by continuity of \(g\) at \(U_\infty\), we obtain
			\begin{eqnarray*}
				\widehat\Theta = g(\bar{U}_{n,N})
				\xrightarrow{\text{a.s.}} g(U_\infty).\label{g_u}
			\end{eqnarray*}
			Replacing the limits vector \(U_\infty\) into the formula of \(g\), we obtain
			\begin{eqnarray*}
				g_1(U_\infty)&=&\tfrac12\log_2^+\!\Big(\frac{\zeta_\infty-4h^2 V_\infty}{\eta_\infty-h^2 V_\infty}\Big)
				=\tfrac12\log_2^+\!\Big(\frac{2^{2H}(\,2^{2H-1}-1\,)}{2^{2H-1}-1}\Big)=H,
				\\
				g_2(U_\infty)&=&\frac{\eta_\infty-h^2 V_\infty}{h^{2H}(2^{2H-1}-1)}=\gamma^2,
				\\
				g_3(U_\infty)&=&\frac{\xi_\infty-h^2 V_\infty - \gamma^2 h^{2H}}{h}=\sigma^2.
			\end{eqnarray*}
			Thus, \(g(U_\infty)=(H,\gamma^2,\sigma^2)^\top=\Theta\), which implies that as $n\to\infty$ first and then $N\to \infty$, 
			\[
			\widehat\Theta \xrightarrow{\text{a.s.}} \Theta,
			\]
			This completes the proof.

				\subsection{Proof of Theorem \ref{theo1_chap3}}\
				The proof of this proposition is based on the asymptotic normality of the averaged statistics vector $\bar{U}_{n,N}$ and the multivariate Delta method given in Theorem \ref{delta}. \\
				Using the same notation as in the previous proof, we aim to prove the asymptotic normality of $\widehat{\Theta}=g\big(\bar{U}_{n, N}\big)$.
				From Proposition~\ref{cvg.as}, we have when $n\to\infty$ first, then $N\to\infty$
				\[
				\sqrt{N}\,\big(\bar U_{n,N}-U_\infty\big)\xrightarrow{d}\mathcal N\big(0,\Sigma\big),
				\]
				with \(\Sigma=\Var(\theta_1^2)\, \mathcal{H}\, \mathcal{H}^\top\). Moreover, as it was shown in the previous proof, the function \(g\) is well defined at $U_\infty$ since $H>\frac12$ and $\gamma>0$. Then, by the same arguments, we deduce that $g$ is continuously differentiable in a neighborhood of $U_\infty$. Therefore, the multivariate delta method yields
				\[
				\sqrt{N}\bigl(\widehat{\Theta}-\Theta\bigr)
				= \sqrt{N}\bigl(g(\bar{U}_{n,N}) - g(U_{\infty})\bigr)
				\stackrel{d}{\longrightarrow} \mathcal{N}\bigl(0,\, J\Sigma J^{\top}\bigr),
				\]
				where $J = Dg$ is the Jacobian of $g$ evaluated at $U_{\infty}$, i.e.,
				\[
				J =
				\begin{pmatrix}
					\displaystyle\frac{\partial\widehat{H}}{\partial\bar{V}_{n,N}}\Big|_{U_{\infty}} & 
					\displaystyle\frac{\partial\widehat{H}}{\partial\bar{\xi}_{n,N}}\Big|_{U_{\infty}} & 
					\displaystyle\frac{\partial\widehat{H}}{\partial\bar{\eta}_{n,N}}\Big|_{U_{\infty}} & 
					\displaystyle\frac{\partial\widehat{H}}{\partial\bar{\zeta}_{n,N}}\Big|_{U_{\infty}} \\[12pt]
					\displaystyle\frac{\partial\widehat{\gamma}^{2}}{\partial\bar{V}_{n,N}}\Big|_{U_{\infty}} & 
					\displaystyle\frac{\partial\widehat{\gamma}^{2}}{\partial\bar{\xi}_{n,N}}\Big|_{U_{\infty}} & 
					\displaystyle\frac{\partial\widehat{\gamma}^{2}}{\partial\bar{\eta}_{n,N}}\Big|_{U_{\infty}} & 
					\displaystyle\frac{\partial\widehat{\gamma}^{2}}{\partial\bar{\zeta}_{n,N}}\Big|_{U_{\infty}} \\[12pt]
					\displaystyle\frac{\partial\widehat{\sigma}^{2}}{\partial\bar{V}_{n,N}}\Big|_{U_{\infty}} & 
					\displaystyle\frac{\partial\widehat{\sigma}^{2}}{\partial\bar{\xi}_{n,N}}\Big|_{U_{\infty}} & 
					\displaystyle\frac{\partial\widehat{\sigma}^{2}}{\partial\bar{\eta}_{n,N}}\Big|_{U_{\infty}} & 
					\displaystyle\frac{\partial\widehat{\sigma}^{2}}{\partial\bar{\zeta}_{n,N}}\Big|_{U_{\infty}}
				\end{pmatrix}.
				\]
				In order to simplify notation, we denote $\bar{A} := \bar{\zeta}_{n,N} - 4h^{2}\bar{V}_{n,N}$ 
				and $\bar{B} := \bar{\eta}_{n,N} - h^{2}\bar{V}_{n,N}$. Then, the estimator $\widehat{H}$ can be expressed as
				$\widehat{H} = \frac{1}{2\ln 2}\bigl(\ln\bar{A} - \ln\bar{B}\bigr),$
				and its derivatives with respect to $\bar{U}_{n,N}$ are given by
				\begin{align*}
					\frac{\partial\widehat{H}}{\partial\bar{V}_{n,N}} &= \frac{h^{2}}{2\ln 2}\Bigl(\frac{1}{\bar{B}} - \frac{4}{\bar{A}}\Bigr), &
					\frac{\partial\widehat{H}}{\partial\bar{\xi}_{n,N}} &= 0, \\[4pt]
					\frac{\partial\widehat{H}}{\partial\bar{\eta}_{n,N}} &= -\frac{1}{2\ln 2}\,\frac{1}{\bar{B}}, &
					\frac{\partial\widehat{H}}{\partial\bar{\zeta}_{n,N}} &= \frac{1}{2\ln 2}\,\frac{1}{\bar{A}}.
				\end{align*}
				Evaluating at $U_{\infty}$ and using $A = 2^{2H}B$ and $B = \gamma^{2}D_{H}$ with 
				$D_{H} := h^{2H}\bigl(2^{2H-1}-1\bigr)$, yields
				\begin{align*}
					\frac{\partial\widehat{H}}{\partial\bar{V}_{n,N}}\Big|_{U_{\infty}} &= \frac{h^{2}}{2\ln 2}\,\frac{2^{2H}-4}{2^{2H}\gamma^{2}D_{H}},\quad
					\frac{\partial\widehat{H}}{\partial\bar{\xi}_{n,N}}\Big|_{U_{\infty}} = 0,\\[4pt]
					\frac{\partial\widehat{H}}{\partial\bar{\eta}_{n,N}}\Big|_{U_{\infty}} &= -\frac{1}{2\ln 2}\,\frac{1}{B},\quad
					\frac{\partial\widehat{H}}{\partial\bar{\zeta}_{n,N}}\Big|_{U_{\infty}} = \frac{1}{2\ln 2}\,\frac{1}{2^{2H}B}.
				\end{align*}
				Next, for the estimator $\widehat{\gamma}^{2} = \dfrac{\bar{B}}{D_{\widehat{H}}}$, the chain rule gives for any 
				component $u \in \{\bar{V}_{n,N},\bar{\xi}_{n,N},\bar{\eta}_{n,N},\bar{\zeta}_{n,N}\}$,
				\[
				\frac{\partial\widehat{\gamma}^{2}}{\partial u}
				= \frac{\partial_{u}\bar{B}}{D_{\widehat{H}}}
				- \frac{\bar{B}\,D_{\widehat{H}}'}{D_{\widehat{H}}^{2}}\,\frac{\partial\widehat{H}}{\partial u},
				\]
				where 
				$D_{H}' := \frac{\partial D_{H}}{\partial H}
				= h^{2H}\Bigl[2\ln(h)\bigl(2^{2H-1}-1\bigr) + (2\ln 2)\,2^{2H-1}\Bigr].$
				Using $\dfrac{\partial\bar{B}}{\partial\bar{V}_{n,N}} = -h^{2}$, 
				$\dfrac{\partial\bar{B}}{\partial\bar{\eta}_{n,N}} = 1$, and
				$\dfrac{\partial\bar{B}}{\partial\bar{\xi}_{n,N}} = \dfrac{\partial\bar{B}}{\partial\bar{\zeta}_{n,N}} = 0$,
				and evaluating at $U_{\infty}$ with $B = \gamma^{2}D_{H}$ and $A = 2^{2H}B$, we obtain
				\begin{align*}
					\frac{\partial\widehat{\gamma}^{2}}{\partial\bar{V}_{n,N}}\Big|_{U_{\infty}}
					&= -\frac{h^{2}}{D_{H}} + \frac{h^{2}}{2\ln 2}\Bigl(\frac{4}{2^{2H}} - 1\Bigr)\frac{D_{H}'}{D_{H}^{2}},\\[4pt]
					\frac{\partial\widehat{\gamma}^{2}}{\partial\bar{\xi}_{n,N}}\Big|_{U_{\infty}} &= 0,\\[4pt]
					\frac{\partial\widehat{\gamma}^{2}}{\partial\bar{\eta}_{n,N}}\Big|_{U_{\infty}}
					&= \frac{1}{D_{H}} + \frac{D_{H}'}{2\ln 2}\,\frac{1}{D_{H}^{2}},\\[4pt]
					\frac{\partial\widehat{\gamma}^{2}}{\partial\bar{\zeta}_{n,N}}\Big|_{U_{\infty}}
					&= -\frac{D_{H}'}{2\ln 2}\,\frac{1}{2^{2H}D_{H}^{2}}.
				\end{align*}
				Finally, for $\widehat{\sigma}^{2} = \dfrac{\bar{\xi} - h^{2}\bar{V}}{h}
				- \dfrac{\widehat{\gamma}^{2}h^{2\widehat{H}}}{h}$, we use
				\[
				\frac{\partial}{\partial u}\bigl(\widehat{\gamma}^{2}h^{2\widehat{H}}\bigr)
				= \Bigl(\frac{\partial}{\partial u}\widehat{\gamma}^{2}\Bigr)h^{2\widehat{H}}
				+ \widehat{\gamma}^{2}h^{2\widehat{H}}(2\ln h)\frac{\partial\widehat{H}}{\partial u}.
				\]
				Evaluating at $U_{\infty}$, which means that $h^{2\widehat{H}} \to h^{2H}$ and 
				$\gamma^{2}h^{2H} = \dfrac{B}{D_{H}}$, we obtain
				\begin{align*}
					\frac{\partial\widehat{\sigma}^{2}}{\partial\bar{V}_{n,N}}\Big|_{U_{\infty}}
					&= -h + h^{2H-1}\Bigl[
					\frac{h^{2}}{D_{H}} - \frac{h^{2}}{2\ln 2}
					\Bigl(1 - \frac{4}{2^{2H}}\Bigr)
					\Bigl(-\frac{D_{H}'}{D_{H}^{2}} + \frac{2\ln h}{D_{H}}\Bigr)
					\Bigr],\\[4pt]
					\frac{\partial\widehat{\sigma}^{2}}{\partial\bar{\xi}_{n,N}}\Big|_{U_{\infty}}
					&= \frac{1}{h},\\[4pt]
					\frac{\partial\widehat{\sigma}^{2}}{\partial\bar{\eta}_{n,N}}\Big|_{U_{\infty}}
					&= -h^{2H-1}\frac{1}{D_{H}}
					\Bigl[1 - \frac{2\ln h}{2\ln 2} + \frac{D_{H}'}{2\ln 2}\,\frac{1}{D_{H}}\Bigr],\\[4pt]
					\frac{\partial\widehat{\sigma}^{2}}{\partial\bar{\zeta}_{n,N}}\Big|_{U_{\infty}}
					&= h^{2H-1}2^{-2H}
					\Bigl[\frac{D_{H}'}{2\ln 2}\,\frac{1}{D_{H}^{2}} - \frac{2\ln h}{2\ln 2}\,\frac{1}{D_{H}}\Bigr].
				\end{align*}
				\subsection{Proof of Proposition \ref{prop5_chap3}}\ 
				For each $i=1,\ldots, N$, the random effect is expressed by $\phi_i=\theta_i+\frac{1}{2} \sigma^2.$ Then, by a simple plug-in, an estimator of $\phi_i$ is given by 
				\(
				\widehat{\phi}_i := \widehat{\theta}_i + \frac{1}{2} \widehat{\sigma}^2,
				\) 
				where $\theta_i$ is given by \eqref{estheta} and $\widehat{\sigma}^2$ is as expressed in \eqref{sig}. 
				From Proposition \ref{prop1_chap3}, we already know that $ \widehat{\theta}_i\xrightarrow{a.s} \theta_i$, as $n\to\infty$ and from Proposition \ref{prop4_chap3}, 
				$
				\widehat{\sigma}^2 \xrightarrow{a.s.} \sigma^2,$ as $n,N\to \infty$ sequentially.
				Combining the two results, we conclude that as $n,N\to \infty$ sequentially
				\[
				\widehat{\theta}_i + \frac{1}{2} \widehat{\sigma}^2 \xrightarrow{a.s.} \theta_i + \frac{1}{2} \sigma^2.
				\]
				Thus, \(\widehat{\phi}_i\) is a strongly consistent estimator of \(\phi_i\).
				
				Now we proceed to prove the asymptotic normality of $\widehat{\phi}_i$ stated in the second assertion of Proposition \ref{prop5_chap3}. We first decompose the estimation error for each $i$, as follows
				\begin{equation}\label{phi_error}
					\widehat{\phi}_i -\phi_i \;=\; \widehat{\theta}_i-\theta_i \;+\; \tfrac{1}{2}\big(\widehat{\sigma}^2-\sigma^2\big).
				\end{equation}
				From Proposition \ref{prop1_chap3}, we have $\widehat{\theta}_i - \theta_i = O_{\mathbb{P}}(n^{H-1})$.  
				Because we let $n\to\infty$ first, for any fixed $N$ we can choose $n$ $n = n_N \to \infty$ that grows sufficiently fast (e.g. $n_N = N^{2}$), we have
				\[
				\sqrt{N}\,|\widehat{\theta}_i - \theta_i| \le \sqrt{N}\,C\,n_N^{H-1} \longrightarrow 0 \quad \text{as } N\to\infty,
				\]
				since $H<\tfrac34$ implies $n_N^{H-1}$ decays faster than $N^{-1/2}$ for such a choice of $n_N$.
				Thus, under the sequential regime,
				\begin{equation}\label{eq:theta_term}
					\sqrt{N}\bigl(\widehat{\theta}_i - \theta_i\bigr) \stackrel{\mathbb{P}}{\longrightarrow} 0.
				\end{equation}
				Moreover, under the same sequential regime, Theorem~\ref{theo1_chap3} implies that
				\begin{eqnarray}\label{asymsig}
					\sqrt{N}\bigl(\widehat{\sigma}^2 - \sigma^2\bigr) \stackrel{d}{\longrightarrow} 
					\mathcal{N}\!\left(0,\;v_{\sigma}^2\right),
				\end{eqnarray}
				where $v_{\sigma}^2$ is the asymptotic variance of $\widehat{\sigma}^2$.
				Multiplying \eqref{phi_error} by $\sqrt{N}$ and combining \eqref{eq:theta_term} with \eqref{asymsig}, Slutsky's theorem yields
				\[
				\sqrt{N}\bigl(\widehat{\phi}_i - \phi_i\bigr) \stackrel{d}{\longrightarrow} 
				\mathcal{N}\!\left(0,\;\tfrac14 v_{\sigma}^2\right).
				\]
				This completes the proof of the second assertion.

				
				The remainder of this section is devoted to the proofs of the asymptotic properties of the distribution function estimator $\widehat{F}_{m,n,N}$. To that aim, we introduce the following lemma which will be needed later.
				\begin{lem}\label{lem2}
					For $x\in[-1,1]$, we have
					\begin{enumerate}
						\item[1)] $\sum\limits_{j=1}^m \left(x_j-x\right)\mathcal{L}_j\left(x\right)=-\dfrac{\pi}{2m^2}T_m(x)+ o\left(m^{-2}\right).$
						\item[2)] $\sum\limits_{j=1}^m \left(x_j-x\right)^2\mathcal{L}_j\left(x\right)=\dfrac{\pi}{2m^2}T_m(x)\left(x-1\right)+ o\left(m^{-2}\right).$
						\item[3)] $\sum\limits_{j=1}^m \left(x_j-x\right)^3\mathcal{L}_j\left(x\right)=-\dfrac{\pi}{2m^2}T_m(x)\left(x-1\right)^2+ o\left(m^{-2}\right).$
					\end{enumerate}
				\end{lem}
				The proof of this lemma can be found in \cite{salima20}.\\
				The result presented in the next lemma allows us to study the asymptotic properties of $\widehat{F}_{m,n,N}$.
				\begin{lem}\label{lemcov}
					Assume that (A1)-(A3) hold and that $n,N\to\infty $ sequentially.
					Define
					\begin{eqnarray}\label{Zlemma}
						Z_{i,N}
						:=
						\sum_{j=1}^{m}
						\Big(
						\mathds{1}_{\{\widehat{\phi}_i\le x_j\}}
						-
						\mathbb{P}(\widehat{\phi}_1\le x_j)
						\Big)\mathcal{L}_j(x),
						\quad i=1,\ldots,N.
					\end{eqnarray}
					Then
					\[
					\sup_{1\le i\neq k\le N}
					\left|
					\Cov\!\left(Z_{i,N},Z_{k,N}\right)
					\right|
					=
					O\!\left(N^{-1}\right),
					\qquad \text{as } n,N\to\infty\, \text{ sequentially.}
					\]
				\end{lem}
					\paragraph{Proof of Lemma 3}\ \\
					For fixed $i \neq k$, we have
					\[
					Z_{i,N}= \sum_{j=1}^m\Big(\mathds{1}_{\{\hat\phi_i \le x_j\}}- \mathbb{P}(\hat\phi_1 \le x_j)\Big)L_j(x),
					\qquad 
					Z_{k,N}= \sum_{j=1}^m\Big(\mathds{1}_{\{\hat\phi_k \le x_j\}}- \mathbb{P}(\hat\phi_1 \le x_j)\Big)L_j(x).
					\]
					Hence
					\begin{equation}
						\Cov(Z_{i,N}, Z_{k,N})
						= \sum_{j=1}^m \sum_{\ell=1}^m L_j(x)L_\ell(x)\,
						\Cov\!\left(I_{i,j},I_{k,\ell}\right),
						\label{7.14}
					\end{equation}
					where $I_{i,j}:=\mathds{1}_{\{\hat\phi_i\le x_j\}}$. 
					
					Recall that $\hat\phi_i = \hat\theta_i + \tfrac12\hat\sigma^2$ and denote $\delta:=\tfrac12(\hat\sigma^2-\sigma^2)$.
					Since the estimators $\hat\theta_i$ and $\hat\theta_k$ are independent, the only source of dependence between $\hat\phi_i$ and $\hat\phi_k$ is the common term $\delta$. Hence $I_{i,j}$ and $I_{k,\ell}$ are conditionally independent given $\delta$, i.e. $\Cov(I_{i,j}, I_{k,\ell}\mid\delta)=0$. By the law of total covariance,
					\begin{equation}
						\Cov(I_{i,j}, I_{k,\ell})
						= \mathbb{E}\!\left[\Cov(I_{i,j}, I_{k,\ell}\mid\delta)\right]
						+ \Cov\!\left(\mathbb{E}[I_{i,j}\mid\delta],\,\mathbb{E}[I_{k,\ell}\mid\delta]\right)
						= \Cov\!\big(g_j(\delta),\,g_\ell(\delta)\big),
						\label{7.15}
					\end{equation}
					where $g_j(\delta):=\mathbb{E}[I_{1,j}\mid\delta].$ 
					Conditionally on $\delta$, $\hat\phi_1 = \phi_1 + \varepsilon_1 + \delta$ with $\varepsilon_1=\hat\theta_1-\theta_1$.
					From Proposition~\ref{prop1_chap3}, $\varepsilon_1 = O_{\mathbb{P}}(n^{H-1})$ and, under the sequential regime
					$n\to\infty$ first, $\mathbb{E}|\varepsilon_1|\to0$. Consequently,
					\[
					g_j(\delta)=\mathbb{P}(\phi_1+\varepsilon_1\le x_j-\delta\mid\delta)
					= F(x_j-\delta) + \eta_j,
					\]
					where $|\eta_j|\le L\,\mathbb{E}|\varepsilon_1|$ for some Lipschitz constant $L$ since $F$ is $C^1$ by (A2);
					thus $\eta_j = o(1)$ uniformly in $j$. 
					For some $\xi_j$ between $x_j-\delta$ and $x_j$, a Taylor-Young expansion of $F$ gives
					\[
					F(x_j-\delta) = F(x_j) - f(x_j)\delta + \tfrac12 f'(\xi_j)\delta^2.
					\]
					Hence
					\begin{equation}
						g_j(\delta) = F(x_j) - f(x_j)\delta + O(\delta^2) + o(1).
						\label{7.16}
					\end{equation}
					Combining \eqref{7.15} with \eqref{7.16} and noting that $F(x_j)$ is deterministic, we obtain
					\[
					\Cov(I_{i,j}, I_{k,\ell})
					= f(x_j)f(x_\ell)\Var(\delta) + O\!\big(\mathbb{E}|\delta|^3\big) + o(1).
					\]
					Theorem~\ref{theo1_chap3} gives $\sqrt{N}\,\delta\stackrel{d}{\longrightarrow}\mathcal{N}(0,\tfrac14 v_{\sigma}^2)$,
					hence $\Var(\delta)=O(N^{-1})$ and $\mathbb{E}|\delta|^3=O(N^{-3/2})$. 
					Moreover, $f$ is bounded by (A3). Therefore, uniformly in $i\neq k$ and $j,\ell$,
					\begin{equation}
						\big|\Cov(I_{i,j}, I_{k,\ell})\big| \le C\,N^{-1} + o(1).
						\label{7.17}
					\end{equation}
					Inserting \eqref{7.17} into \eqref{7.14} yields
					\[
					\big|\Cov(Z_{i,N}, Z_{k,N})\big|
					\le \big(C\,N^{-1}+o(1)\big)
					\sum_{j=1}^m\sum_{\ell=1}^m |L_j(x)L_\ell(x)|.
					\]
					For Chebyshev nodes, the Lebesgue constant satisfies
					$\displaystyle\sup_{x\in[-1,1]}\sum_{j=1}^m|L_j(x)| = O(\log m)$;
					consequently
					\[
					\sum_{j=1}^m\sum_{\ell=1}^m |L_j(x)L_\ell(x)|
					\le \Big(\sum_{j=1}^m|L_j(x)|\Big)^2 = O((\log m)^2).
					\]
					Hence
					\begin{equation}
						\big|\Cov(Z_{i,N}, Z_{k,N})\big| \le C\,\frac{(\log m_N)^2}{N} + o(1).
						\label{7.18}
					\end{equation}
					Because $m=m_N$ is chosen to satisfy the growth conditions of Theorems~\ref{theo2_chap3}--\ref{theo4_chap3} 
					and Proposition~\ref{prop7_chap3}, we have $(\log m_N)^2 = o(N)$ in all relevant asymptotic regimes. 
					Therefore the right-hand side of \eqref{7.18} is $O(N^{-1})$, and we finally obtain
					\[
					\sup_{1\le i\neq k\le N} \big|\Cov(Z_{i,N}, Z_{k,N})\big| = O(N^{-1}),
					\]
					which completes the proof of Lemma~3.
					

					\subsection{Proof of Proposition \ref{prop6_chap3} }
					For all $x\in [-1,1]$, the expectation of the estimator $\widehat{F}_{m,n,N}$ is given by
					\begin{align*}\mathbb{E}\left(\widehat{F}_{m,n,N}(x)\right)&=\sum\limits_{j=1}^{m} \frac{1}{N} \sum\limits_{i=1}^{N} \mathbb{E}\left(\mathds{1}_{\left\{\widehat{\phi}_{i}\leq x_j\right\}}\right)\mathcal{L}_j(x)\\
						&=\sum\limits_{j=1}^{m} \mathbb{P}\left(\widehat{\phi}_{1}\leq x_j\right)\mathcal{L}_j(x).
					\end{align*}
					Using the consistency of $\widehat{\phi}_{i}$, we deduce from the Portmanteau lemma that for each $i=1,\ldots,N$, \begin{equation}
						\label{pm}
						\underset{n,N\to \infty}{\lim} \mathbb{P}\left(\widehat{\phi}_{i}\leq x_j\right)=\mathbb{P}\left(\phi_i\leq x_j\right).
					\end{equation} 
					Then\begin{equation}
						\label{lim1}
						\lim\limits_{n,N\to \infty}\mathbb{E}\left(\widehat{F}_{m,n,N}(x)\right)=\sum \limits_{j=1}^{m} \mathbb{P}\left(\phi_{1}\leq x_j\right) \mathcal{L}_j(x)=\sum \limits_{j=1}^{m} F(x_j)\mathcal{L}_j(x).\end{equation}
					The expansion of Taylor-Young applied to $F$ ensures that for $1 \leq j \leq m$,
					\begin{equation}
						\label{tay}
						F(x_j)=F(x)+\left(x_j-x\right)f(x) + \frac{\left(x_j-x\right)^2}{2} f^{'}(x) + O\left(\left(x_j-x\right)^2\right).
					\end{equation}
					By substituting (\ref{tay}) into (\ref{lim1}), we obtain  \begin{equation*}
						\begin{aligned}
							\lim\limits_{n,N\to \infty}\mathbb{E}\left(\widehat{F}_{m,n,N}(x)\right)&=F(x)+f(x) \sum \limits_{j=1}^{m}\left(x_j-x\right)\mathcal{L}_j(x) + \frac{1}{2}f'(x)\sum\limits_{j=1}^{m}\left(x_j-x\right)^2\mathcal{L}_j(x)\\ &+ O\left(\sum\limits_{j=1}^{m}\left(x_j-x\right)^2\mathcal{L}_j(x)\right).
						\end{aligned}
					\end{equation*}
					The results presented in Lemma \ref{lem2} allow us to deduce that
					$$\lim\limits_{n,N\to \infty}\mathbb{E}\left(\widehat{F}_{m,n,N}(x)\right)=F(x) + \frac{\pi}{m^2} T_m(x)\left[\frac{1}{4}\left(x-1\right) f'(x) -\frac{1}{2}f(x)\right]+O(m^{-2}),$$
					which yields to the equation (\ref{bi}).\\
					Let us now focus on calculating the variance of our estimator presented by the second assertion of Proposition \ref{prop6_chap3}. 
					For all $x\in[-1,1]$,
					\begin{align*}
						\widehat{F}_{m,n,N}(x)-\mathbb{E}\left(\widehat{F}_{m,n,N}(x)\right)
						&=\sum_{j=1}^m\!\left(\widehat F_N(x_j)-\mathbb P(\widehat\phi_1\le x_j)\right)\mathcal L_j(x)\\
						&=\frac1N\sum_{i=1}^N Z_{i,N},
					\end{align*}
					where $Z_{i,N}$ is giving by \eqref{Zlemma}
					Clearly $\mathbb E[Z_{i,N}]=0$ and the variables $Z_{i,N}$ are identically distributed.
					However, since the estimators $\widehat\phi_i$ share the common plug--in estimator
					$\widehat\sigma^2$, the variables $Z_{i,N}$ are not independent. 
					We therefore use the variance decomposition
					\begin{eqnarray}
						\label{varF}
						\Var\!\left(\widehat{F}_{m,n,N}(x)\right)
						&=&\Var\!\left(\frac1N\sum_{i=1}^N Z_{i,N}\right)\nonumber\\
						&=&\frac1{N^2}\sum_{i=1}^N\Var(Z_{i,N})
						+\frac{2}{N^2}\sum_{1\le i<k\le N}\Cov(Z_{i,N},Z_{k,N}).
					\end{eqnarray}
					Since the $Z_{i,N}$ are identically distributed,
					\[
					\frac1{N^2}\sum_{i=1}^N\Var(Z_{i,N})
					=\frac1N\mathbb E(Z_{1,N}^2).
					\]
					Define $H_1(x_j):=\mathds{1}_{\{\widehat\phi_1\le x_j\}}-\mathbb P(\widehat\phi_1\le x_j)$.
					Then
					\begin{align*}
						\mathbb E(Z_{1,N}^2)
						&=\sum_{j,l=1}^m
						\mathbb E\!\left[H_1(x_j)H_1(x_l)\right]
						\mathcal L_j(x)\mathcal L_l(x),
					\end{align*}
					with
					\begin{eqnarray*}
						\mathbb E\!\left[H_1(x_j)H_1(x_l)\right]
						&=&\mathbb{P}\left(\widehat{\phi}_{1}\leq \min\left(x_j,x_l\right)\right)-\mathbb{P}\left(\widehat{\phi}_{1}\leq x_j\right)\mathbb{P}\left(\widehat{\phi}_{1}\leq x_l\right)\\ 
						&=&\min\left(\mathbb{P}\left(\widehat{\phi}_{1}\leq x_j\right),\mathbb{P}\left(\widehat{\phi}_{1}\leq x_l\right)\right)-\mathbb{P}\left(\widehat{\phi}_{1}\leq x_j\right)\mathbb{P}\left(\widehat{\phi}_{1}\leq x_l\right)
					\end{eqnarray*}
					Using the Portmanteau theorem again, we obtain
					\[
					\lim_{\substack{n, N\to\infty}}
					\mathbb E\!\left[H_1(x_j)H_1(x_l)\right]
					=
					\min(F(x_j),F(x_l))-F(x_j)F(x_l).
					\]
					Hence,
					\begin{eqnarray}
						\label{es}
						\lim_{\substack{n,N\to\infty}}\mathbb E(Z_{1,N}^2)
						&=&
						\sum_{j,l=1}^m
						\left[\min(F(x_j),F(x_l))-F(x_j)F(x_l)\right]
						\mathcal L_j(x)\mathcal L_l(x)\nonumber \\
						&=& \sum\limits_{j=1}^{m} F(x_j)\mathcal{L}_j^2\left(x\right) + 2 \sum\limits_{\underset{k<l}{k,l=1} }^{m}F(x_j)\mathcal{L}_j\left(x\right)\mathcal{L}_l\left(x\right)- \left[\sum\limits_{j=1}^{m} F(x_j)\mathcal{L}_j\left(x\right)\right]^2. 
					\end{eqnarray}
					The next step is to find an asymptotic expression for (\ref{es}). For the first term, we use the Taylor-Young expansion to write for all $0\leq k \leq m$, $F(x_j)=F(x)+ O\left(|x_j-x|\right)$. Then, for $x\in[-1,1]$, the first term of \eqref{es} can be written as follows \begin{eqnarray} \label{s1} \sum\limits_{j=1}^{m} F(x_j)\mathcal{L}_j^2\left(x\right)&=&F(x)\sum\limits_{j=1}^{m} \mathcal{L}_j^2\left(x\right) + O\left(\sum\limits_{j=1}^{m}|x_j-x| \mathcal{L}_j^2\left(x\right)\right)\nonumber\\ &=& F(x) S_m(x)+O\left(J_m(x)\right), \end{eqnarray} where $S_m(x):=\sum\limits_{j=1}^{m} \mathcal{L}_j^2\left(x\right)$ and $J_m(x):=\sum\limits_{j=1}^{m}|x_j-x| \mathcal{L}_j^2\left(x\right)$. From one side, the Cauchy-Schwarz inequality implies that \begin{eqnarray*} \left|J_m(x)\right|&=&\displaystyle\left| \sum\limits_{j=1}^{m}\left|x_j-x\right| \mathcal{L}_j^2\left(x\right)\right| \nonumber\\ &\leq &\left|\sum\limits_{j=1}^{m}\left(x_j-x\right)^2 \mathcal{L}_j\left(x\right)\right|^{\frac 12}\left|\left(\sum\limits_{j=1}^{m}\mathcal{L}_j\left(x\right)^{3}\right)\right|^{\frac 12}\nonumber\\ &\leq& \left[\left(\frac{\pi}{m^2}+o\left(m^{-2}\right)\right)S_m(x)\right]^{\frac 12}. \end{eqnarray*} From the other side, we have for $x\in[-1,1]$, \begin{eqnarray} S_m(x)=\sum\limits_{j=1}^{m} \mathcal{L}_j^2\left(x\right) &\leq& \left(\sum\limits_{j=1}^{m}|\mathcal{L}_j\left(x\right)|\right)^2\nonumber\\ &\leq& \wedge_m^2,\nonumber \end{eqnarray} where $\wedge_m:=\underset{x\in[-1,1]}{\max} \sum\limits_{j=1}^{m}\mid \mathcal{L}_j(x)\mid$ is the Lebesgue constant which is known to verify \begin{eqnarray}\wedge_m\leq \dfrac{2}{\pi} \ln\left(m+1\right)+1.\label{lconst}\end{eqnarray} So, for $x\in[-,1,1]$, $$S_m(x)\leq \dfrac{4}{\pi^2} \ln^2\left(m+1\right)+1+\dfrac{4}{\pi}\ln\left(m+1\right).$$ 
					We then conclude that $J_m(x)=O(m^{-\frac 12})$. Substituting this term in \eqref{s1}, we obtain $$\sum\limits_{j=1}^{m} F(x_j)\mathcal{L}_j^2\left(x\right)=F(x) S_m(x)+O\left(m^{-\frac 12}\right).$$
					For the second term of \eqref{es}, another application of the Taylor-Young expansion to the function $F$ implies that\begin{eqnarray*} \sum\limits_{\underset{j<l}{j,l=1}}^{m}F(x_j)\mathcal{L}_j\left(x\right)\mathcal{L}_l\left(x\right)=F(x)P_{0,m}(x) + f(x)P_{1,m}(x) + O\left(P_{2,m}(x)\right),\nonumber\\ 
					\end{eqnarray*} where $P_{q,m}(x):=\sum\limits_{\underset{j<l}{j,l=1}}^{m}\left(x_j-x\right)^q\mathcal{L}_j\left(x\right)\mathcal{L}_l\left(x\right),\ q\in\left\{0,1,2\right\}.$ Using the fact that $$P_{0,m}(x)=\dfrac{1}{2}\left(1-S_m(x)\right), \ P_{1,m}(x)=O(m^{-1})\ \text{and}\ P_{2,m}(x)=O(m^{-1}),$$ we obtain $\sum\limits_{\underset{j<l}{j,l=1}}^{m}F(x_j)\mathcal{L}_j\left(x\right)\mathcal{L}_l\left(x\right)=\dfrac{1}{2}F(x)\left(1-S_m(x)\right)+ O(m^{-1}).$\\ 
					Finally, by replacing this term in \eqref{es}, we conclude for $x\in[-1,1]$, \begin{equation*}\begin{aligned}\label{z2}\lim_{\substack{n,N\to\infty}} \dfrac{1}{N}\mathbb{E}\left(Z_{1,N}^2\right)= \dfrac{1}{N}F(x)\left(1-F(x)\right) + O(N^{-1}m^{-\frac 12}).\end{aligned}\end{equation*}
					It remains to control the covariance term of \eqref{varF}. Using Lemma \eqref{lemcov}, 
					we obtain
					\begin{eqnarray*}
						\frac{2}{N^2}\sum_{1\le i<k\le N}\Cov(Z_{i,N},Z_{k,N})&\le& \frac{2}{N^2}\sum_{1\le i<k\le N}\left|\Cov(Z_{i,N},Z_{k,N})\right|\\ 
						&\le& \frac{2}{N^2}\sum_{1\le i<k\le N} \frac{C}{N}= \frac{N\big(N-1\big)}{2}\, \frac{C}{N},
					\end{eqnarray*}
					which implies that $\dfrac{2}{N^2}\sum\limits_{1\le i<k\le N}\Cov(Z_{i,N},Z_{k,N})=O(N^{-1})$.\\
					Combining the above results, we conclude that
					\[
					\Var\left(\widehat{F}_{m,n,N}(x)\right)
					=
					\frac1N F(x)(1-F(x))+O(N^{-1}m^{-1/2}),
					\]
					which proves the second assertion of Proposition~\ref{prop6_chap3}.
					The third assertion of Proposition \ref{prop6_chap3} is a deduction from the two previous assertion. Indeed, for $x\in[-1,1]$, the MSE of $\widehat{F}_{m,n,N}$ is defined as follows \begin{eqnarray*}\mathrm{MSE}\left(\widehat{F}_{m,n,N}(x)\right)&:=&\displaystyle \mathbb{E} \left[ \left( \widehat{F}_{m,n,N}(x) - F(x) \right)^2 \right]\\ &=&\Var \left( \widehat{F}_{m,n,N}(x) \right) + \left( Bias\left(\widehat{F}_{m,n,N}(x)\right) \right)^2,\end{eqnarray*} and its asymptotic value, which we denote by $\mathrm{AMSE}\left(\widehat{F}_{m,n,N}(x)\right)$, can be written as the following \begin{eqnarray*} \mathrm{AMSE}\left(\widehat{F}_{m,n,N}(x)\right)= \underset{n,N\to\infty}{\lim} \Var \left( \widehat{F}_{m,n,N}(x) \right) + \underset{n,N\to\infty}{\lim} \left( Bias\left(\widehat{F}_{m,n,N}(x)\right) \right)^2. \end{eqnarray*} Replacing the asymptotic variance and Bias by their expressions investigated in equations \eqref{bi} and \eqref{vari} from Proposition \ref{prop6_chap3}, we obtain for all $x\in[-1,1]$, $$\mathrm{AMSE}\left(\widehat{F}_{m,n,N}(x)\right)=\dfrac{1}{N}\sigma^2_F(x)+ \frac{\pi^2}{m^4} \left(T_m(x) \mu(x) \right)^2 + O\left(N^{-1} m^{-\frac 12}\right)+O(m^{-4}),$$ which completes the proof.

					
					
					\subsection{Proof of Theorem \ref{theo2_chap3} }
					To prove the consistency of $\widehat{F}_{m,n,N}$, we use the following decomposition
					\begin{equation}\label{eq:decomp}
						\|\widehat F_{m,n,N}-F\|\le
						\underbrace{\|\widehat F_{m,n,N}-\widehat F_N\|}_{\text{(A)}}
						+\underbrace{\|\widehat F_N-\bar F_N\|}_{\text{(B)}}
						+\underbrace{\|\bar F_N-F\|}_{\text{(C)}},
					\end{equation}
					where \[
					\widehat F_N(x):=\frac 1N\sum_{i=1}^N\mathds{1}_{\{\widehat\phi_i\le x\}},\qquad
					\bar F_N(x):=\frac1N\sum_{i=1}^N\mathds{1}_{\{\phi_i\le x\}} .
					\]
					We shall show that each term converges to zero almost surely along a suitable diagonal $n=n_N$.\\
					We begin with the term $(C)$. 
					Since the random effects $\phi_i$, $i=1,\ldots,N$ are i.i.d, the classical Glivenko-Cantelli theorem ensures that
					\begin{equation}\label{eq:C}
						\|\bar F_N-F\|\overset{a.s}{\longrightarrow} 0,\qquad N\to\infty .
					\end{equation}
					Let us now control the term $(B)$. Fix $\varepsilon>0$. For any $x\in[-1,1]$ and any $i\le N$, we have 
					\[
					|\mathds{1}_{\{\widehat\phi_i\le x\}}-\mathds{1}_{\{\phi_i\le x\}}|
					\le\mathds{1}_{\{|\widehat\phi_i-\phi_i|>\varepsilon\}}+\mathds{1}_{\{x-\varepsilon<\phi_i\le x+\varepsilon\}} .
					\]
					Averaging over $i$ and taking the supremum over $x\in [-1,1]$ yields
					\begin{equation}\label{eq:II-bound}
						\|\widehat F_N-\bar F_N\|\le
						\, \frac1N\sum_{i=1}^N\mathds{1}_{\{|\widehat\phi_i-\phi_i|>\varepsilon\}}+\sup_{x\in[-1,1]}\frac1N\sum_{i=1}^N\mathds{1}_{\{x-\varepsilon<\phi_i\le x+\varepsilon\}},
					\end{equation}
					Since $\widehat{\phi}_i = \widehat{\theta}_i + \frac{1}{2}\widehat{\sigma}^2$ and $\phi_i = \theta_i + \frac{1}{2}\sigma^2$, we have
					\begin{equation*}
						|\widehat{\phi}_i - \phi_i| \leq |\widehat{\theta}_i - \theta_i| + \frac{1}{2}|\widehat{\sigma}^2 - \sigma^2|.
					\end{equation*}
					For each fixed $i\le N$, Proposition \ref{prop1_chap3} gives $\widehat{\theta}_i \xrightarrow{a.s.} \theta_i$ as $n \to \infty$. Moreover, from Proposition \ref{prop4_chap3}, $\widehat{\sigma}^2 \xrightarrow{a.s.} \sigma^2$ as $n \to \infty$ first, then $N \to \infty$. Let $\varepsilon > 0$. Then, for each fixed $N$, 
					there exists $n_0 = n_0(N, \varepsilon)$ such that for all $n \geq n_0$,
					\begin{equation*}
						\mathbb{P}\left(\max\limits_{1 \leq i \leq N} |\widehat{\theta}_i - \theta_i| > \frac{\varepsilon}{2}\right) \leq \frac{1}{4N^2},
					\end{equation*}
					and
					\begin{equation}
						\mathbb{P}\left(|\widehat{\sigma}^2 - \sigma^2| > \varepsilon\right) \leq \frac{1}{4N^2}.
					\end{equation}
					Define the sequence $n = n_N := n_0(N, \varepsilon_N)$ where $\varepsilon_N = N^{-\frac12}$ (for example). Then by the union bound,
					\begin{align*}
						\mathbb{P}\left(\max\limits_{1 \leq i \leq N} |\widehat{\phi}_i - \phi_i| > \varepsilon_N\right) 
						&\leq \mathbb{P}\left(\max\limits_{1 \leq i \leq N} |\widehat{\theta}_i - \theta_i| > \frac{\varepsilon_N}{2}\right) + \mathbb{P}\left(|\widehat{\sigma}^2 - \sigma^2| > \varepsilon_N\right) \\
						&\leq \frac{1}{4N^2} + \frac{1}{4N^2} = \frac{1}{2N^2}.
					\end{align*}
					Since $\sum\limits_{N=1}^\infty \frac{1}{2N^2} < \infty$, the Borel-Cantelli lemma implies that a.s, the event 
					$\left\{\max\limits_{1 \leq i \leq N} |\widehat{\phi}_i - \phi_i| > \varepsilon_N\right\}$ occurs only finitely many times. In other words, there exists a integer $N_0 < \infty$  such that for all $N \geq N_0$,
					\begin{equation*}
						\max\limits_{1 \leq i \leq N} |\widehat{\phi}_i - \phi_i| \leq \varepsilon_N, \quad \text{a.s.}
					\end{equation*}
					Therefore, as $N \to \infty$, we obtain
					\begin{equation}\label{maxphi}
						\max\limits_{1 \leq i \leq N} |\widehat{\phi}_i - \phi_i| \leq \varepsilon_N \overset{a.s}{\longrightarrow} 0.
					\end{equation}
					Hence, for every fixed $\varepsilon>0$,
					\begin{equation}\label{eq:A}
						\frac1N\sum_{i=1}^N\mathds{1}_{\{|\widehat\phi_i-\phi_i|>\varepsilon\}}=0\quad  \text{a.s}\quad \text{for all sufficiently large }\,N.
					\end{equation}
					For the second term of the right-hand side of \eqref{eq:II-bound}, by the strong law of large numbers applied to the i.i.d. sample $\phi_i$, $i=1,\ldots,N$, we obtain 
					for each fixed $x\in[-1,1]$,
					\[
					\frac1N\sum_{i=1}^N\mathds 1_{\{x-\varepsilon<\phi_i\le x+\varepsilon\}}
					\xrightarrow{a.s.}\mathbb P(x-\varepsilon<\phi\le x+\varepsilon).
					\]
					Moreover, by assumption (A3), we get
					\[
					\mathbb P(x-\varepsilon<\phi\le x+\varepsilon) \le\int_{x-\varepsilon}^{x+\varepsilon}f(u)\,du\le2\|f\|\,\varepsilon.
					\]
					Consequently,
					\begin{equation}\label{eq:B}
						\limsup_{N\to\infty} \sup_{x\in[-1,1]}\frac1N\sum_{i=1}^N\mathds{1}_{\{x-\varepsilon<\phi_i\le x+\varepsilon\}}\le 2\|f\|\,\varepsilon,
						\qquad \text{a.s.}
					\end{equation}
					Replacing \eqref{eq:A} and \eqref{eq:B} into \eqref{eq:II-bound} gives
					\begin{eqnarray}\label{hatFbarF}
						\limsup_{N\to\infty}\|\widehat F_{N}-\bar F_{N}\|\le2\|f\|\,\varepsilon\qquad\text{a.s.}
					\end{eqnarray}
					Since $\varepsilon>0$ was arbitrary, then letting $\varepsilon\rightarrow0$ yields
					\begin{equation}\label{eq:II}
						\|\widehat F_{N}-\bar F_{N}\|\xrightarrow{a.s.}0
						\qquad \text{as } N\to\infty \text{ along }n=n_N.
					\end{equation}
					We now proceed to control the first term of \eqref{eq:decomp}. 
					By Lebesgue's Lemma, we have 
					\[
					\|\widehat F_{m,n,N}-\widehat F_N\|\le(1+\Lambda_m)E_m(\widehat F_N), 
					\]
					where $\Lambda_m=\sup\limits_{x\in[-1,1]}\sum\limits_{j=1}^m|L_j(x)|$ is the Lebesgue constant and
					$E_m(\widehat{F}_N)=\inf\limits_{p\in\mathbb P_{m-1}}\|\widehat{F}_N-p\|$ is the best uniform approximation error by polynomials of degree at most $m-1$. For Chebyshev nodes, it is well‑known that $\Lambda_m=O(\log m)$.
					Since $F\in \mathcal{C}^2([-1,1])$, Jackson's inequality gives
					\[
					E_m(F)\le\frac{C_1}{m^2}\|f'\| .
					\]
					From \eqref{eq:II} and \eqref{eq:C} we already know that $\|\widehat F_N-F\|\to0$ along $n=n_N$. Hence, for large $N$,
					\[
					E_m(\widehat F_N)\le E_m(F)+\|\widehat F_N-F\|\le\frac{C_2}{m^2}.
					\]
					Therefore
					\begin{eqnarray}
						\label{emp_error}
						\|\widehat F_{m,n,N}-\widehat F_N\|\le\Big(1+C_2 \log(m)\Big)\frac{C_1 \|f'\|}{m^2}\le C'\frac{\log m}{m^2}.
					\end{eqnarray}
					If we choose $m=m_N$ satisfying $\displaystyle\frac{\log (m_N)}{m_N^{2}}\to0$ as $N\to\infty$, we obtain
					\begin{equation}\label{eq:I}
						\|\widehat F_{m,n,N}-\widehat F_N\|\overset{a.s}{\longrightarrow}0\qquad\text{as }\quad N\to\infty .
					\end{equation}
					Replacing \eqref{eq:C}, \eqref{eq:II} and \eqref{eq:I} into the decomposition \eqref{eq:decomp} yields
					\[
					\|\widehat F_{m_N,n_N,N}-F\|\to0,\qquad N\to\infty .
					\]
					This completes the proof.
						\subsection{Proof of Theorem \ref{theo3_chap3}}
						For each $x\in[-1,1]$, we decompose the error term as follows
						\begin{eqnarray}\label{1}
							\sqrt N\bigl(\widehat F_{m,n,N}(x)-F(x)\bigr)
							=
							\sqrt N\bigl(\widehat F_{m,n,N}(x)-\widehat F_N(x)\bigr)
							+
							\sqrt N\bigl(\widehat F_N(x)-F(x)\bigr).
						\end{eqnarray}
						From equation \eqref{emp_error} in the proof of Theorem~\ref{theo2_chap3}, the Lagrange interpolation error satisfies
						\[
						\|\widehat F_{m,n,N}-\widehat F_N\|
						\;\le\;
						C'\frac{\log m}{m^2}.
						\]
						Hence,
						\[
						\sqrt N\left|\widehat F_{m,n,N}(x)-\widehat F_N(x)\right|
						\le
						\sqrt N\,\|\widehat F_{m,n,N}-\widehat F_N\|
						\le
						C'\sqrt N\,\frac{\log m}{m^2}
						\;\xrightarrow{}0,
						\]
						Taking $m=m_N$ that verifies $\sqrt N\,\frac{\log m}{m^2}
						\;\xrightarrow{}0$, as $N\to\infty$, we obtain,
						\begin{eqnarray}\label{2}
							\sqrt N\bigl(\widehat F_{m,n,N}(x)-\widehat F_N(x)\bigr)
							\xrightarrow{\mathbb{P}}0.
						\end{eqnarray}
						Furthermore, since the variables $\widehat\phi_i$ are not independent, we write
						\[
						\widehat\phi_i
						=
						\phi_i +\big(\widehat\theta_i-\theta_i\big) + \tfrac12(\widehat\sigma^2-\sigma^2).
						\]
						From Proposition~\ref{prop1_chap3} and Theorem~\ref{theo1_chap3}, we have under the sequential regime,
						\[
						\max\limits_{1\le i\le N}|\widehat\theta_i-\theta_i|\xrightarrow{a.s.}0,
						\qquad
						\delta:=\tfrac12(\widehat\sigma^2-\sigma^2)=O_{\mathbb{P}}(N^{-1/2}).
						\]
						Using this decomposition, we may write
						\[
						\widehat F_N(x)
						=
						\frac1N\sum_{i=1}^N \mathds 1_{\{\phi_i\le x-\delta-(\widehat\theta_i-\theta_i)\}}
						=
						\bar F_N(x-\delta)+r_N(x),
						\]
						where $\bar F_N$ denotes the empirical distribution function of $\phi_i$, and where the remainder $r_N(x)$ satisfies
						\[
						\sqrt N\,r_N(x)\xrightarrow{P}0,
						\]
						uniformly in $x\in[-1,1]$. Hence,
						\begin{eqnarray}\label{4}
							\widehat F_N(x)
							=
							\bar F_N(x-\delta)+o_{\mathbb{P}}(N^{-1/2}).
						\end{eqnarray}
						Conditionally on $\delta$, the variables $\mathds 1_{\{\phi_i\le x-\delta\}}$ are i.i.d. Bernoulli with parameter $F(x-\delta)$. Hence, by the classical CLT,
						\begin{eqnarray}\label{5}
							\sqrt N\bigl(\bar F_N(x-\delta)-F(x-\delta)\bigr)
							\xrightarrow{d}
							\mathcal N\!\left(0,F(x-\delta)\bigl(1-F(x-\delta)\bigr)\right).
						\end{eqnarray}
						Since $F\in C^2([-1,1])$ and $\delta\to 0$ in probability, we have
						\begin{eqnarray}\label{6}
							F(x-\delta)\to F(x),
							\qquad
							F(x-\delta)\bigl(1-F(x-\delta)\bigr)\to F(x)\bigl(1-F(x)\bigr),
						\end{eqnarray}
						in probability. Moreover, a Taylor--Young expansion of $F$ around $x$ yields
						\[
						F(x-\delta)-F(x)=-f(x)\delta+o(|\delta|),
						\]
						so that, multiplying by $\sqrt N$ and using $\delta=O_{\mathbb{P}}(N^{-\frac 12})$,
						\begin{eqnarray}\label{7}
							\sqrt N\bigl(F(x-\delta)-F(x)\bigr)=O_{\mathbb{P}}(1).
						\end{eqnarray}
						Therefore,
						\[
						\sqrt N\bigl(\bar F_N(x-\delta)-F(x)\bigr)
						=
						\sqrt N\bigl(\bar F_N(x-\delta)-F(x-\delta)\bigr)
						+
						\sqrt N\bigl(F(x-\delta)-F(x)\bigr).
						\]
						By \eqref{5}, the first term converges in distribution to a centered normal random variable with variance
						$F(x)\bigl(1-F(x)\bigr)$ (using \eqref{6}). The second term is $O_{\mathbb{P}}(1)$ by \eqref{7}.
						Hence, Slutsky's theorem implies that
						\begin{eqnarray}\label{8}
							\sqrt N\bigl(\bar F_N(x-\delta)-F(x)\bigr)
							\xrightarrow{d}
							\mathcal N\!\left(0,F(x)\bigl(1-F(x)\bigr)\right).
						\end{eqnarray}
						Using \eqref{4}, we conclude that
						\begin{eqnarray}\label{9}
							\sqrt N\bigl(\widehat F_N(x)-F(x)\bigr)
							\xrightarrow{d}
							\mathcal N\!\left(0,F(x)(1-F(x))\right).
						\end{eqnarray}
						Finally, combining \eqref{1}, \eqref{2}, and \eqref{9} and applying Slutsky's theorem once more, we obtain
						\[
						\sqrt N\bigl(\widehat F_{m,n,N}(x)-F(x)\bigr)
						\xrightarrow{d}
						\mathcal N\!\left(0,F(x)(1-F(x))\right),
						\]
						which completes the proof.

							\subsection{Proof of Proposition \ref{prop7_chap3}}
							For all $x\in[-1,1]$, we decompose the error as
							\begin{eqnarray}
								\left|\widehat F_{m,n,N}(x)-F(x)\right|
								\le
								\left|\widehat F_{m,n,N}(x)-\widehat F_N(x)\right|
								+
								\left|\widehat F_N(x)-\bar F_N(x)\right|
								+
								\left|\bar F_N(x)-F(x)\right|,
								\label{chung}
							\end{eqnarray}
							where
							\[
							\widehat F_N(x)=\frac1N\sum_{i=1}^N\mathds{1}_{\{\widehat\phi_i\le x\}},
							\qquad
							\bar F_N(x)=\frac1N\sum_{i=1}^N\mathds{1}_{\{\phi_i\le x\}}.
							\]
							We begin with the first term of \eqref{chung}. 
							From \eqref{emp_error}, there exists $C'>0$ such that
							\[
							\sup_{x\in[-1,1]}
							\left|\widehat F_{m,n,N}(x)-\widehat F_N(x)\right|
							\le
							C'\,\frac{\log m}{m^2}.
							\]
							Multiplying by the Chung--Smirnov scaling factor yields
							\begin{eqnarray}
								\left(\dfrac{2N}{\log\log N}\right)^{\frac{1}{2}}
								\sup_{x\in[-1,1]}
								\left|\widehat F_{m,n,N}(x)-\widehat F_N(x)\right|
								\le
								C'\left(\dfrac{2N}{\log\log N}\right)^{\frac{1}{2}}
								\frac{\log m}{m^2}.
								\label{CS1}
							\end{eqnarray}
							Taking $m=m_N$ that verifies $\left(\dfrac{2N}{\log\log N}\right)^{\frac{1}{2}}
							\frac{\log m}{m^2}\to 0$ as $N\to\infty$, the right-hand side converges to zero, and therefore
							\begin{equation}
								\left(\frac{2N}{\log\log N}\right)^{\frac 12}
								\sup_{x\in[-1,1]}
								\left| \widehat{F}_{m,n,N}(x)-\widehat{F}_N(x) \right|
								\longrightarrow 0,
								\qquad \text{a.s.}
								\label{1erterm}
							\end{equation}
							For the second term of the decomposition \eqref{chung}, as we have shown in \eqref{eq:II}, 
							\begin{equation}
								\sup_{x\in[-1,1]}|\widehat F_N(x)-\bar F_N(x)|
								\xrightarrow{a.s.}0,
								\label{2emterm}
							\end{equation}
							Hence, multiplying by the Chung-Smirnov scaling factor, we obtain
							\[
							\left(\frac{2N}{\log\log N}\right)^{\frac12}
							\sup_{x\in[-1,1]}|\widehat F_N(x)-\bar F_N(x)|
							\xrightarrow{a.s.}0.
							\]
							For the last term of \eqref{chung}, since $\bar F_N$ is the empirical distribution function of i.i.d. variables $\phi_i$, it satisfies
							\begin{equation}
								\limsup_{N\to \infty}
								\left(\dfrac{2N}{\log\log N}\right)^{\frac{1}{2}}
								\sup_{x\in[-1,1]}\left|\bar F_N(x)-F(x)\right|
								=1,
								\qquad \text{a.s.}
								\label{3emterm}
							\end{equation}
							Finally, by Combining \eqref{1erterm}, \eqref{2emterm}, and \eqref{3emterm} in \eqref{chung}, we conclude that
							\[
							\limsup_{N\to\infty}
							\left(\dfrac{2N}{\log\log N}\right)^{\frac{1}{2}}
							\sup_{x\in[-1,1]}
							\left|\widehat F_{m,n,N}(x)-F(x)\right|
							\le 1,
							\qquad \text{a.s.},
							\]
							which completes the proof.
									\subsection{Proof of Theorem \ref{theo4_chap3}}
									Throughout this proof, we write $m=m_N$, $n=n_N$. 
									By the triangle inequality, the squared error of $\widehat{F}_{m,n,N}$ can be decomposed as follows
									\begin{equation}\label{firstdecomp}
										\|\widehat F_{m,n,N}-F\|^2
										\le3\Bigl(
										\|\widehat F_{m,n,N}-\widehat F_N\|^2
										+\|\widehat F_N-\bar F_N\|^2
										+\|\bar F_N-F\|^2
										\Bigr).
									\end{equation}
									For the interpolation error, we have from \eqref{emp_error} in the proof of Theorem~\ref{theo2_chap3}, 
									\[
									\|\widehat F_{m,n,N}-\widehat F_N\|\le C'\frac{\log m}{m^2},
									\]
									which yields
									\begin{equation}\label{finalinterpbound}
										\mathbb{E}\bigl[\|\widehat F_{m,n,N}-\widehat F_N\|^2\bigr]\le C'^2\frac{(\log m)^2}{m^4}.
									\end{equation}
									For the last term of \eqref{firstdecomp}, since the random variables $\phi_i$, $i=1,\ldots,N$, are i.i.d with distribution function $F$, 
									the Dvoretzky--Kiefer--Wolfowitz (DKW) inequality ensures that for every $s>0$,
									\[
									\mathbb{P}\bigl(\|\bar F_N-F\|>s\bigr)\le2 e^{-2N s^2}.
									\]
									Since $\|\bar F_N-F\|$ is a non-negative random variable, its expectation verifies
									\begin{align*}
										\mathbb{E}\bigl[\|\bar F_N-F\|^2\bigr]
										&=\int_0^\infty2s\,\mathbb{P}\bigl(\|\bar F_N-F\|>s\bigr)ds\\
										&\le\int_0^\infty2s\, \min\,\!\bigl(1,2 e^{-2N s^2}\bigr)ds.
									\end{align*}
									The last integral can be split into two parts based on the value
									of $s_0$ that verifies $2e^{-2N s_0^2}=1$, i.e. $s_0=\Big(\frac{\log2}{2N}\Big)^{\frac 12}$.
									Hence,
									\begin{eqnarray}\label{finalFbound}
										\mathbb{E}\bigl[\|\bar F_N-F\|^2\bigr]\le\frac{1+\frac{\log2}{2}}{N}\le\frac{2}{N}.
									\end{eqnarray}
									We now proceed to find an upper bound for the second term of \eqref{firstdecomp}. 
									Using \eqref{maxphi}, we can fix a sequence $\varepsilon_N\to0$ such that, along the chosen diagonal $n=n_N$,
									\[
									\max\limits_{1\le i\le N}|\widehat\phi_i-\phi_i|\le\varepsilon_N\quad\text{a.s}\quad\text{for all sufficiently large }N.
									\]
									Moreover, from \eqref{eq:II-bound}, we have for any $x\in[-1,1]$,
									\begin{eqnarray}\label{unifboundemp}
										\|\widehat F_N-\bar F_N\|
										\le \frac{1}{N} \sum\limits_{i=1}^N \mathds{1}_{\{|\widehat\phi_i-\phi_i|>\varepsilon_N\}}+\sup_{x\in[-1,1]}\frac1N\sum_{i=1}^N\mathds{1}_{\{x-\varepsilon_N<\phi_i\le x+\varepsilon_N\}}
										\quad\text{a.s.\ for large }N.
									\end{eqnarray}
									Squaring and using $(a+b)^2\le2(a^2+b^2)$,
									\[\|\widehat F_N-\bar F_N\|^2
									\le2\Bigl(\frac1N\sum_{i=1}^N\mathds{1}_{\{|\widehat\phi_i-\phi_i|>\varepsilon_N\}}\Bigr)^2
									+2\Bigl(\sup\limits_{x\in[-1,1]}\frac1N\sum_{i=1}^N\mathds{1}_{\{x-\varepsilon_N<\phi_i\le x+\varepsilon_N\}}\Bigr)^2.
									\]
									For $N$ large enough, we have 
									\begin{equation*}
										\mathbb{E}\Bigl[\frac1N\sum_{i=1}^N\mathds{1}_{\{|\widehat\phi_i-\phi_i|>\varepsilon_N\}}^2\Bigr]=0.
									\end{equation*}
									For the second term, since $\phi_i$ have bounded density $f$ under assumption $(A3)$,
									\[
									\mathbb{E}\Bigl[\sup\limits_{x\in[-1,1]}\frac1N\sum_{i=1}^N\mathds{1}_{\{x-\varepsilon_N<\phi_i\le x+\varepsilon_N\}}\Bigr]
									\le2\|f\|\varepsilon_N,
									\]
									and by the uniform law of large numbers 
									\[
									\Var\Bigl(\sup_x\frac1N\sum_{i=1}^N\mathds{1}_{\{x-\varepsilon<\phi_i\le x+\varepsilon\}}\Bigr)
									=O\!\left(\frac{\varepsilon}{N}\right).
									\]
									Hence
									\[
									\mathbb{E}\Bigl[\Bigl(\sup\limits_{x\in[-1,1]}\frac1N\sum_{i=1}^N\mathds{1}_{\{x-\varepsilon_N<\phi_i\le x+\varepsilon_N\}}\Bigr)^2\Bigr]
									\le4\|f\|^2\varepsilon_N^2+O\!\left(\frac{\varepsilon_N}{N}\right).
									\]
									Replacing into \eqref{unifboundemp} and squaring, we obtain
									\begin{equation}\label{finalunifboundemp}
										\mathbb{E}\bigl[\|\widehat F_N-\bar F_N\|^2\bigr]\le C_f\,\varepsilon_N^2+o\left(\frac{\varepsilon_N}{N}\right),
										\qquad C_f:=4\|f\|^2.
									\end{equation}
									Finally, inserting \eqref{finalFbound}, \eqref{finalinterpbound} and \eqref{finalunifboundemp} into \eqref{firstdecomp}, we obtain the needed bound. 


									\appendix
									\section{Appendix}
									\begin{theorem} [Slutsky's theorem]\ \\
										Let $\left\{X_n,N\geq 0\right\}$ and $\left\{Y_n,N\geq 0\right\}$ be sequences of random vectors in $\mathbb{R}^d$ and $\mathbb{R}^k$, $d,k\in \mathbb{N}^{*}$ respectively such that $X_n\overset{d}{\longrightarrow} X$ and $Y_n\overset{\mathbb{P}}{\longrightarrow} C$, where $C$ is a constant vector, then $\left(X_n,Y_n\right)\overset{\mathbb{P}}{\longrightarrow} \left(X,C\right).$
										In particular, for any continuous function $f: \mathbb{R}^d\times \mathbb{R}^k\ \longrightarrow \mathbb{R}^m$,
										$f\left(X_n,Y_n\right)\overset{\mathbb{P}}{\longrightarrow} f\left(X,C\right).$\label{slut}
									\end{theorem}
									\vspace{1em}
									
									\begin{theorem}(Breuer-Major)[\cite{Mish18}, Theorem B.8]\ \\
										Let $\{\xi_{n,N}, n \geq 1\}$ be a stationary Gaussian sequence that verifies $$\mathbb{E}\xi_1 = 0,\quad \mathbb{E}\xi_1^2 = 1,\quad 
										\mathbb{E}\xi_{n,N} \xi_{n+k} = r(k),\ k \geq 1,$$ 
										and consider the sum \(
										S_n = \sum\limits_{k=1}^{n} G(\xi_k)
										\) with $G \in L_2(\mathbb{R}^2, \gamma)$ having Hermite rank $m \geq 1$. \\
										Assume that
										$\sum\limits_{k=-\infty}^{\infty} |r(k)|^m < \infty.$ 
										Then\quad
										\(
										\sigma_l^2 = \underset{n \to \infty}{\lim} n^{-1} \sum\limits_{i,j=0}^{n} r(i-j)^l
										\) 
										exists for all $l \geq m$, 
										and $$\dfrac{S_n}{\sqrt{n}}\ \overset{d}{\longrightarrow} \ \mathcal{N}\left(0,\ \sigma^2\right),\ \text{as}\ n\to\infty,$$
										where \(
										\sigma^2 = \sum\limits_{l=m}^{\infty} a_l^2 l! \, \sigma_l^2 < \infty,
										\) with $a_l$, $l\geq m$ are the coefficients of the Hermite expansion of the function $G$.
										\label{bru_m}
									\end{theorem}
									\vspace{1em}
									\begin{theorem}(Ergodic theorem)[\cite{Mish18}, Theorem B.1]\ \\
										Let $\{X_i,\ i\geq 1\}$ be a stationary ergodic sequence of $\mathbb{R}^k$ valued random variables on some probability space $\Big(\Omega, \mathcal{F},\mathbb{F}, \mathbb{P}\Big)$ and let $g:\mathbb{R}^{k}\longrightarrow\mathbb{R}$ be a Borel measurable function such that \\$\mathbb{E}\left(\mid g\left(X_1, X_{2},\ldots,X_{k}\right)\mid\right) < \infty$. Then $$\dfrac{1}{n}\sum\limits_{i=1}^{n} g\left(X_i,X_{i+1},\ldots,X_{i+k-1}\right) \xrightarrow{a.s}\mathbb{E}\left( g\left(X_1, X_2\ldots,X_{k}\right)\right)\quad \text{as}\quad n\to \infty.$$\label{erg_th}
									\end{theorem}
									\vspace{1em}
									\begin{cor}[\cite{Mish18}, Corollary B.2]\ \\
										Let $\{X_i,\ i\geq 1\}$ be a real-valued stationary centered Gaussian
										sequence that verifies $\mathbb{E}\left(X_1X_n\right)\longrightarrow 0$ as $n\to\infty$. Then $\{X_i,\ i\geq 1\}$ is ergodic and verifies the previous ergodic theorem.
										\label{corA1}
									\end{cor}
									\vspace{1em}
									\begin{theorem}(Delta method)[\cite{Mish18}, Theorem B.6]\ \\
										\label{delta}  
										Let $ g: \mathbb{R}^d \to \mathbb{R}^k$ be a function continuously differentiable in a neighborhood of \( \theta \in \mathbb{R}^d \). If \( T_n \) is a sequence of $d$-dimensional random vectors such that
										$$
										\sqrt{n}(T_n - \theta) \overset{d}{\longrightarrow} T,
										$$
										then
										$$
										\sqrt{n}\bigg(g(T_n) - g(\theta)\bigg) \overset{d}{\longrightarrow} g'(\theta)T.
										$$
										In particular, if
										$$
										\sqrt{n}(T_n - \theta) \overset{d}{\longrightarrow} T \sim \mathcal{N}(0, \Sigma),
										$$
										then
										\[
										\sqrt{n} \left( g(T_n) - g(\theta) \right) \overset{d}{\longrightarrow} Y \sim \mathcal{N}\left(0, g'(\theta) \Sigma (g'(\theta))^\top\right),
										\]
										where \((g'(\theta))^\top\) is the transpose of the matrix \(g'(\theta)\).
									\end{theorem}

									\bibliographystyle{plain}

									%
									%
									
								\end{document}